\documentclass[journal,twoside,web]{ieeecolor}
\usepackage{generic}
\usepackage{amsmath,amssymb,amsfonts}
\usepackage{algorithmic}
\usepackage{graphicx}
\usepackage{subfigure}
\usepackage{amsmath,bm}
\usepackage{textcomp}
\usepackage{multirow}
\usepackage{tikz}
\graphicspath{{images/}}
\usepackage{amsmath}
\usepackage{latexsym, amssymb}
\usepackage{amsbsy}
\usepackage{amsopn, amstext}
\usepackage{ifpdf,hyperref}
\usepackage{cancel, color}
\usepackage{epstopdf}
\def\BibTeX{{\rm B\kern-.05em{\sc i\kern-.025em b}\kern-.08em
		T\kern-.1667em\lower.7ex\hbox{E}\kern-.125emX}}
\def\J{{\bf 1}}

\usepackage[]{algorithm2e}
\usepackage{tikz}
\graphicspath{{images/}}
\usepackage{amsmath}
\usepackage{latexsym, amssymb}
\usepackage{graphicx}
\usepackage{amsmath, amsbsy}
\usepackage{amsopn, amstext}
\usepackage{ifpdf,hyperref}
\usetikzlibrary{calc,trees,positioning,arrows,chains,shapes.geometric,decorations.pathreplacing,decorations.pathmorphing,shapes,matrix,shapes.symbols}

\usepackage{tikz}
\usepackage{amsmath}
\usepackage{bm}
\usepackage{latexsym, amssymb, amsmath, amsbsy,amsopn, amstext}
\usepackage{graphicx,hyperref,booktabs}
\usepackage{cancel, color}
\usepackage{epstopdf}
\graphicspath{{images/}}

\DeclareMathOperator{\Span}{Span}
\DeclareMathOperator{\Col}{Col}

\DeclareMathOperator{\lcm}{lcm}

\def\cal{\mathcal}

\def\pa{\partial}

\def\ra{\rightarrow}

\def\a{\alpha}
\def\b{\beta}
\def\d{\delta}

\def\D{\Delta}

\def\0{{\bf 0}}

\def\J{{\bf 1}}

\newcommand{\R}{{\mathbb R}}

\def\dsum{\mathop{\sum}\limits}

\newtheorem{thm}{Theorem}[section]
\newtheorem{cor}[thm]{Corollary}
\newtheorem{dfn}[thm]{Definition}
\newtheorem{prp}[thm]{Proposition}
\newtheorem{exa}[thm]{Example}
\newtheorem{alg}[thm]{Algorithm}
\newtheorem{rem}[thm]{Remark}
\newtheorem{lem}[thm]{Lemma}

\begin{document}

\title{Universal Solution to Kronecker Product Decomposition}
\author{Daizhan Cheng
	\thanks{D. Cheng is with the Key Laboratory of Systems and Control, Academy of Mathematics and Systems Science, Chinese Academy of Sciences, Beijing 100190, P.R.China, e-mail: dcheng@iss.ac.cn}
	\thanks{This work is supported partly by the National Natural Science Foundation of China (NSFC) under Grant 62350037.}
}

\maketitle

\begin{abstract}

This paper provides a general solution for the Kronecker product decomposition (KPD) of vectors, matrices, and hypermatrices. First, an algorithm, namely, monic decomposition algorithm (MDA), is reviewed. It consists of a set of projections from a higher dimension Euclidian space to its factor-dimension  subspaces. It is then proved that the KPD of vectors is solvable, if and only if, the project mappings provide the required decomposed vectors. Hence it provides an easily verifiable necessary and sufficient condition for the KPD of vectors. Then an algorithm is proposed to calculate the least square error approximated decomposition.
Using it finite times a finite sum (precise) KPD of any vectors can be obtained. Then the swap matrix is used to make the elements of a matrix re-arranging, and then provides a method to convert the KPD of matrices to its corresponding KPD of vectors. It is proved that the KPD of a matrix is solvable, if and only if, the KPD of its corresponding vector is solvable. In this way, the necessary and sufficient condition, and the followed algorithms for approximate and finite sum KPDs for matrices are also obtained. Finally, the permutation matrix is introduced and used to convert the KPD of any hypermatrix to KPD of its corresponding vector. Similarly  to matrix case, the necessary and sufficient conditions for solvability and the techniques for vectors and matrices are also applicable for hypermatrices, though some additional algorithms are necessary. Several numerical examples are included to demonstrate this universal KPD solving method.

\end{abstract}


\IEEEpeerreviewmaketitle

\section{Introduction}

Consider a matrix $A\in {\cal M}_{m\times n}$, where $m=m_1m_2$ and $n=n_1n_2$. Then KPD problem consists of several types.
\begin{itemize}
\item[(i)] Exact KPD\cite{van93}

Find   $B\in {\cal M}_{m_1\times n_1}$ and $C\in {\cal M}_{m_2\times n_2}$ such that
\begin{align}\label{1.1}
A=B\otimes C.
\end{align}

\item[(ii)] Approximated KPD\cite{van93}

When the precise decomposition does not exist, find the nearest approximated decomposition.
\begin{align}\label{1.101}
A\approx B\otimes C.
\end{align}

\item[(iii)] Finite sum KPD\cite{cai23}

Find  $B_i\in {\cal M}_{m_1\times n_1}$ and $C_i\in {\cal M}_{m_2\times n_2}$, $i\in [1,k]$, such that
\begin{align}\label{1.2}
A=\dsum_{i=1}^k \lambda_i B_i\otimes C_i,\quad \lambda_i\in \R.
\end{align}

\end{itemize}

In recent years, the KPD of matrices has attracted a considerable attention from different societies, such as computer  and numerical analysis, systems science, signal processing, etc., because  its wide applications, covering (a) various problems from numerical linear algebra \cite{cha00} and the references therein; (b) analysis and synthesis of logical dynamic systems \cite{wei23}; (c) analysis of medical imaging data\cite{fen24}; (d) direction-of-arrival estimation in traffic and manufacturing systems \cite{wan21,joe24}; (e)room acoustic impulse responses \cite{dog22};(f) system identification \cite{pal18}, just to mention a few.

Particularly, the KPD technique has also been found many applications in artificial intelligence, because it can
greatly reducing the total number of parameters of the model. This is particularly important for large scale models.  The related researches include (a)
    training block-wise sparse matrices for machine learning \cite{zhu24}; (b) GPT compression \cite{eda21,ayo24}; (c) CNN compression \cite{ham22}; (c) image process \cite{wu23b}, etc.

Although many theoretical and numerical solutions for KPD have been obtained for various situations and to meet different demands, but in general it is still a challenging problem. For instance, to our best knowledge, theoretically an easily verifiable necessary and sufficient condition for KPD solvability is still open, and numerically a universal efficient computation is still developing.

Although there has been some efforts in KPD of hypermatrices\cite{bat16}, the main efforts have been focused on the KPD of matrices. Strangely enough, the KPD of vectors are not discussed enough.

Recently, an easily computable method, called the monic decomposition algorithm (MDA) has been proposed by us  \cite{chepr}  to solve the  decomposition of vectors. The decomposable vectors are called the hypervectors. MDA consists of a set of mappings, which map a vector to its factor-dimensional subspaces. When a vector is decomposable, the projected vectors form its factor components.

MDA provides not only a necessary and sufficient condition for the MDA of vectors, but also a constructive algorithm to realize the factorization. It is obvious that the algorithm is robust, because it consists of a set of linear mappings, which are continuous. Hence if a higher dimensional vector is near decomposable, the MDA will  provide an approximate decomposition.

If a vector is not decomposable, this paper develops an algorithm, which use the approximate decomposition obtained from MDA as its initial values and  provides a decomposition with least square error. Secondly, if the finite sum decomposition is required, we can use the algorithm on the error vector (the difference between original vector and its approximated hypervector) to get second term approximate decomposition. Keep doing finite times, the algorithm yields the finite sum (precise) KPD of the original vector.
In one word, the paper develops a numerical method to solve the three kinds of KPD for vectors.

Next, as a key step, we prove a result which shows that the KPD of a matrix is equivalent to the KPD of its corresponding vector.
Hence the necessary and sufficient condition for the decomposability, the exact KPD, the least square error approximate KPD, and the precise finite sum KPD of matrices are also solved completely.

Finally, we consider the KPD of hypermatrices. For a hypermatrix of order $d$, we suggest a standard type of matrix expression as its default matrix form. Then the Kronecker product of hypermatrices can be performed via their matrix forms.
When $d=3$, the hypermatrices are also called the cubic matrices. The a default matrix expression may be considered as the one used  in papers of t-product\cite{kil13,arz18}.

Using permutation matrix, the elements of a hypermatrix can be re-arranged. By the light of the permutation, the KPD of hypermatrices can also be reduced to the KPD of its corresponding vectors. That means {\bf the necessary and sufficient condition for the decomposability, the exact KPD, the least square error approximate KPD, and the precise finite sum KPD of an arbitrary hypermatrix  are also solved completely}.

The rest of this paper is organized as follows: Section 2 contains some preliminaries including the semi-tensor product (STP) of matrices, hypervectors, algebraic structure of index set, and hypermatrix with its matrix expression. The permutation matrix is also introduced. Section 3 shows the difference between matrix set and matrix form. It is emphasized that a hypermatrix is not a matrix but a matrix set. Section 4 considers the KPD of vectors. The  MDA is introduced. Then the exact KPD is discussed. Using MDA, a necessary and sufficient condition for exact KPD is obtained, which is constructive via a decomposition algorithm. The approximate MDA and finite sum (precise) MDA of vectors are then following. The kernel technique for KPD is developed in this section. The MDA of matrices is discussed in Section 5. The key point of this KPD is to convert it to the corresponding KPD of vectors, and the equivalence is proved. Then the necessary and sufficient condition, the exact, approximate, and the finite sum KPD of matrices are all obtained.  In Section 6 the formal definition of Kronecker product of hypermatrices  is defined, and three kinds of Kronecker products are presented, which are outer product, partition-based product, and paired product. The KPD of the first two products are solved in this secton. Section 7 devotes to the paired product of hypermatrices and its KPD. The normal matrix form of order $d$ hypermatrices is constructed. The Kronecker product of hypermatrices is described via the matrix forms of the factor hypmatrices. Then the KPD of an arbitrary order $d$ hypmatrix is solved via converting it to the KPD of its corresponding vector form. It follows that the necessary and sufficient condition, the exact, approximate, and the finite sum KPD of hypermatrices are all solvable. Section 8 is a brief conclusion.

Before ending this section a list of notations is attached as follows.

\begin{enumerate}
\item $\R^n$: $n$ dimensional real Euclidean space.
\item ${\cal M}_{m\times n}$: the set of $m\times n$ matrices.
\item $\Col_i(A)$: i-th column of $A$.
\item $\lcm(a,b)$: least common multiple of $a$ and $b$.
\item $A\ltimes B$: semi-tensor product of $A$ and $B$.

\item  $\R^{n_1\times n_2\times \cdots \times n_d}$: the $d$-th order hypermatrices of dimensions $n_1\times\cdots\times n_d$.
\item  $\R^{n_1\ltimes n_2\ltimes \cdots \ltimes n_d}$: the $d$-th order hypervectors of dimensions $n_1,n_2,\cdots,n_d$.
\item ${\cal A}, {\cal B}$, etc.: hypermatrices.

\item $\vec{i}=Id(i;n)$: index $i\in [0,n]$,

\item $M^{{\bf j}\times {\bf k}}({\cal A})$: matrix expression of ${\cal A}$, which has rows labeled by ${\bf j}$ and columns labeled by ${\bf k}$.

\item $V_r(A)$ ($V_c(A)$): row stacking form (column stacking form) of $A$.

\item $[x]$: Integer part of $x$.

\item $[a,b]$: the set of integers $a\leq i \leq b$.
\item $\d_n^i$: the $i$-th column of the identity matrix $I_n$.
%
%
\item $\D_n:=\left\{\d_n^i\;|\; i=1,\cdots,n\right\}$.
\item $\d_n[i_1,\cdots,i_s]:=\left[\d_n^{i_1},\cdots,\d_n^{i_s}\right]$.
\item $\J_{\ell}:=(\underbrace{1,1,\cdots,1}_{\ell})^\mathrm{T}$.
%
%
\item ${\bf S}_d$: $d$-th order permutation group.

%
\item $\otimes$:  Kronecker product of matrices.
\item $\ltimes$: semi-tensor product (STP) of matrices.
%
%
%
\end{enumerate}

\section{Preliminaries}

This section provides some necessary mathematical preliminaries concerning the STP, hypermatrices and hypervectors.

\subsection{STP of Matrices}

\begin{dfn}\label{d2.1.1}  \cite{che12} Let $A\in {\cal M}_{m\times n}$, $B\in {\cal M}_{p\times q}$, and $t=\lcm(n,p)$. The STP of $A$ and $B$ is defined as
\begin{align}\label{2.1.1}
		A\ltimes B=\left(A\otimes I_{t/n}\right) \left(B\otimes I_{t/p}\right).
	\end{align}
\end{dfn}

The STP is a generalization of the classical matrix product, i.e., when $n=p$,  $A\ltimes B=AB$. Because of this, we  omit the symbol $\ltimes$ in most cases. Some fundamental properties of STP are sumarized as follows.

\begin{prp}\label{p2.1.2}
	\begin{itemize}
		\item[(i)] (Associativity) Let $A,B,C$ be three matrices of arbitrary dimensions. Then $(A\ltimes B)\ltimes C=A\ltimes (B\ltimes C)$.
		\item[(ii)] (Distributivity) Let $A,B$ be two matrices of same dimension, $C$ is of arbitrary dimension. Then $(A + B)\ltimes C=A\ltimes C + B\ltimes C$, $C\ltimes (A + B)=C\ltimes A + C\ltimes B$.
		\item[(iii)] Let $A,B$ be two matrices of arbitrary dimensions. Then $(A \ltimes B)^T=B^T\ltimes A^T$. If $A$, $B$ are invertible, then $(A \ltimes B)^{-1}=B^{-1}\ltimes A^{-1}$.
\item[(iv)] If $x,~y$ are two column vectors, then $x\ltimes y=x\otimes y$; if $x,~y$ are two row vectors,
then $x\ltimes y=y\otimes x$.
	\end{itemize}
\end{prp}

Through this paper the default matrix product is assumed to be STP, and the product symbol $\ltimes$ is mostly omitted.

\subsection{Hypervectors}

\begin{dfn}\label{d2.2.1}
\begin{itemize}
\item[(i)] Let $x\in \R^n$, where $n=\prod_{i=1}^rn_i$, if there exist $x_i\in \R^{n_i}$, $i\in [1,r]$ such that
$x=\ltimes_{i=1}^rx_i$, then $x$ is called a hypervector with order $r$ and dimension $n_1\times \cdots\times n_r$.  $x_i$, $i\in [1,r]$ are components of $x$.
\item[(ii)] The set of hypervectors of order $r$ and  dimension $n_1\times \cdots\times n_r$ is denoted by $\R^{n_1\ltimes \cdots\ltimes n_r}$.
\item[(iii)] Assume $0\neq x\in\R^{n_1\ltimes \cdots\ltimes n_r}$. Let
$e=\min_{s}\{x_s\neq 0\}$, called the head index of $x$, and $x_e$ is called the head value of $x$.  $x$ is called a monic vector, if $x_e=1$.
\end{itemize}
\end{dfn}

\begin{prp}\label{p2.2.2}
Let $x\in \R^t$ be a column vector, then for a matrix $A$ with arbitrary dimension,
\begin{align}\label{2.2.1}
x\ltimes A=(I_t\otimes A)\ltimes x.
\end{align}
\end{prp}

\begin{dfn}\label{d2.2.3} The swap matrix (with dimension $m\times n$) is defined by
\begin{align}\label{2.2.2}
W_{[m,n]}:=\left[I_n\d_m^1,I_n\d_m^2,\cdots,I_n\d_m^m \right]\in {\cal M}_{mn\times mn}.
\end{align}
\end{dfn}

\begin{prp}\label{p2.2.4} Let $x\in \R^m$ and $y\in \R^n$ be two column vectors. Then
\begin{align}\label{2.2.3}
W_{[m,n]}xy=yx.
\end{align}
\end{prp}

\begin{lem}\label{l2.2.5}\cite{chepr} Assume $0\neq x\in\R^{n_1\ltimes \cdots\ltimes n_r}$ with the head index $e$, and head value $x_e=a\neq 0$. Then there exists a unique decomposition
\begin{align}\label{2.2.4}
x=a\ltimes_{i=1}^rx_i,
\end{align}
where $x_i$, $i\in[1,r]$ are all monic.
\end{lem}

Using Lemma \ref{l2.2.5}, the following result is obvious.

 \begin{prp}\label{p2.2.6} Let $x=\ltimes_{i=1}^dx_i$, where $x_i\in \R^n$, $i\in [1,d]$. Then there exists a unique subspace
 \begin{align}\label{2.2.5}
S(x):=\Span\{x_1,\cdots, x_d\}\subset \R^n.
\end{align}
 \end{prp}

\subsection{Index Monoid}

We first review the index which we used to label  the elements in a hypermatrix, and giving the set of indexes an algebraic structure.

\begin{dfn}\label{dkd.1.1.1}
\begin{itemize}
\item[(i)] Let $i\in [1,m]$ be an index, where $m>0$ is called the index length. Denote it by
\begin{align}\label{kd.1.1.1}
Id(i;m).
\end{align}
When the index length is clear or ignored,  $Id(i;n)$ is briefly denoted by $\vec{i}$.
\item[(ii)] Let  $i\in [1,m]$ and  $j\in [1,n]$. The product of this two indexes is denoted by
\begin{align}\label{kd.1.1.2}
Id(k;mn)=Id(i;m)\times Id(j;n),
\end{align}
or briefly $\vec{k}=\vec{i}\vec{j}$, where
\begin{align}\label{kd.1.1.3}
k=(i-1)*n+j,\quad i\in [1,m],~j\in [1,n].
\end{align}
\end{itemize}
\end{dfn}

The physical meaning of the product index is that if a set of data is labeled by double index $(i,j)$, where $i\in \vec{i}=Id(i;m)$ and $j\in \vec{j}=Id(j;n)$ as
$D=\{a_{i,j}\;|\; i\in [1,m]. j\in [1,n]\}$, then the elements are arranged in the order of
their product $\vec{k}=Id(k;mn)$. That is, in the single index we have
$$
a_{i,j}=a^k,\quad k=(i-1)*n+j.
$$
To emphasize this relationship, we briefly say that
$$
Id(k;mn)=Id(i;m)\times Id(j;n):=Id(i,j;m,n),
$$
or briefly $\vec{k}=\vec{i}\vec{j}$.

It is easy to verify the following result by definition.

\begin{prp}\label{pkd.1.1.2} The product of indexes is associative. That is.
\begin{align}\label{kd.1.1.4}
(\vec{i}\vec{j})\vec{k}=\vec{i}(\vec{j}\vec{k}).
\end{align}
\end{prp}

Consider ${\cal A}=(a_{i_1,\cdots,i_d})\in \R^{n_1\times \cdots\times n_d}$.
The elements of ${\cal A}$, denoted by
$$
D_{{\cal A}}:=\left\{a_{i_1,\cdots,i_d}\;|\; i_k\in [1,n_k],\; k\in [1,d]\right\},
$$
are arranged in the alphabetic order. In other words, the elements are labeled by
$$
Id(i_1,\cdots,i_d;n_1,\cdots,n_d)=Id(i_1;n_1)\times \cdots \times Id(i_d;n_d).
$$
or briefly, labeled by $\vec{i}=\vec{i}_1\cdots\vec{i}_d$.

 Consider the vector form of ${\cal A}$. Let ${\bf n}=\prod_{i=1}^dn_i$. Using the product (single) index
\begin{align}\label{kd.1.1.6}
Id(j;{\bf n})=\prod_{s=1}^d Id(i_s;n_s)
\end{align}
to label the elements, we have
$$
V({\cal A})=(a^1,a^2,\cdots,a^{\bf n})^{\mathrm{T}}.
$$
Then the conversion relationship can be obtained by a straightforward computation as follows.

\begin{prp}\label{pkd.1.1.3} Consider the index convention (\ref{kd.1.1.6}).
Then
\begin{itemize}
\item[(i)] $\vec{i}\ra \vec{j}$:
\begin{align}\label{kd.1.1.7}
\begin{array}{l}
j=(i_1-1)n_2n_3\cdots n_d+(i_2-1)n_3\cdots n_d\\
~~+\cdots+(i_{d-1}-1)n_d+i_d\\
~~=(\cdot((i_1-1)n_2+(i_2-1))n_3+\cdots\\
~~+(i_{d-1}-1)n_d+i_d.\\
\end{array}
\end{align}
\item[(ii)] $\vec{j}\ra \vec{ i}$:

\begin{align}\label{kd.1.1.8}
\begin{cases}
j_0:=j-1,\\
i_k=j_{d-k}-\left[\frac{j_{d-k}}{n_{k}}\right]+1,\\
j_{d-k+1}=\left[\frac{j_{d-k}}{n_{k}}\right],\quad k=d,d-1,\cdots,2,\\
i_1=j_{d-1}+1.\\
\end{cases}
\end{align}
Where $[x]$ is the integer part of $x$.
\end{itemize}f
\end{prp}

If $\vec{d}=Id(d;1)$, then it becomes a dummy index. We denote the dummy index by $\emptyset$. It follows that
for any index set $\emptyset$ satisfies
$$
\vec{i}\emptyset=\emptyset \vec{i}=\vec{i}.
$$
Then the $\emptyset$ becomes the identity under the index product.  Hence we have

\begin{prp}\label{pkd.1.1.4} The set of indexes with their product is a monoid (i.e., a semi-group with identity).
\end{prp}

\begin{exa}\label{ekd.1.1.4} Consider ${\cal A}=(a_{i_1,\cdots,i_3})\in \R^{3\times 2\times 5}$.
\begin{itemize}
\item[(i)]
Given $V({\cal A})=(a^1,\cdots,a^{{\bf n}})^{\mathrm{T}}$
where ${\bf n}=3*2*5=30$.

(a) Consider $a_{2,1,4}$:
$$
j=(2-1)*10+(1-1)*5+4=14.
$$
That is $a_{2,1,4}=a^{14}$.

(b) Consider $a^{20}$:
$$
\begin{array}{l}
j_0=20-1=19,\\
i_3=19-[19/5]*5+1=5,\\
j_1=[19/5]=3,\\
i_2=3-[3/2]*2+1=2,\\
j_2=1,\\
i_1=j_2+1=2.
\end{array}
$$
Then we have $a^{20}=a_{2,2,5}$.

\item[(ii)]
Given $M^{\{i_1,i_3\}\times \{i_2\}}({\cal A})=(a^{r,s})\in {\cal M}_{15\times 2}$.

(a) Consider $a_{2,1,4}$:
$$
\begin{array}{l}
s=i_2=1,\\
r=(i_1-1)*5+i_3=9,
\end{array}
$$
Hence, $a_{2,1,4}=a^{9,1}$.

(b) Consider $a^{11,2}$. Then $i_2=s=2$. To use formula (\ref{kd.1.1.3} we now consider $i_1$ as the $i_1$ in the formula, and $i_3$ as the $i_2$ in the formula, then we have
$$
\begin{array}{l}
j_0=r-1=10,\\
i_3=j_0-[j_0/5]*5+1=1,\\
j_1=2,\\
i_1=j_1+1=3.
\end{array}
$$
 Hence, $a^{11,2}=a_{3,2,1}$.
\end{itemize}
\end{exa}
\subsection{Hypermatrices}

\begin{dfn}\label{d2.3.1}
\begin{itemize}
\item[(i)] Let $\vec{i}=Id(i;m)$ and  $\vec{j}=Id(j,n)$. Then the product $\vec{k}=\vec{i}\vec{j}$, denoted by $\vec{k}=Id(i,j;m.n)$, stands for the fact that the double index $(i,j)$ is equivalent to the single index $k$. Then
\begin{align}\label{2.3.0}
\begin{array}{l}
k=(i-1)n+j,\\
i=\left[\frac{k-1}{n}\right],~\mbox{and}~j=k-in.\\
\end{array}
\end{align}
\item[(ii)] The multi-index $\vec{k}=\vec{i}_1\vec{i}_2\cdots \vec{i_d}$, denoted by $
Id(i_1,\cdots,i_d;n_1,\cdots,n_d)$, is inductively defined by
$$
\begin{array}{l}
Id(i_1,\cdots,i_k;n_1,\cdots,n_k)=Id(i_1,\cdots,i_{k-1};\\
n_1,\cdots,n_{k-1})Id(i_k;n_k).
\end{array}
$$
\item[(iii)] A set of data
\begin{align}\label{2.3.1}
{\cal A}:=\left\{a_{i_1,\cdots.i_d}\;|\; i_k\in [1,n_k],\;k\in [1,d]\right\}
\end{align}
is said to be labeled by $\vec{i}=Id(i_1,\cdots,i_d;n_1,\cdots,n)_d)$,  which means the elements of ${\cal A}$ are arranged in the alphabetic order as
$$
a_{1,\cdots,1.1},\cdots,a_{1,\cdots,1,n_d},a_{1,\cdots,2,1},\cdots, a_{n_1,\cdots,n_{d-1},n_d}.
$$

\item[(iv)] A set of multi-labeled data as the ${\cal A}$ in (\ref{2.3.1}) is called a hypermatrix of order $d$ and dimension $n_1\times \cdots\times n_d$.

\item[(vi)] The set of hypermatrices with order $d$ and  dimension $n_1\times \cdots\times n_d$ is denoted by
$$
\R^{n_1\times\cdots\times n_d}.
$$
\end{itemize}
\end{dfn}

Note that the definition of hypermatrix here is essentially the same as the one in \cite{lim13}.

To apply the methods and results developed in matrix theory to hypermatrices, the following  technique is fundamental, which converts a hypermatrix into a matrix form.
\begin{dfn}\label{d2.3.2}
\begin{itemize}
\item[(i)] Let ${\cal A}=(a_{i_1,\cdots,i_d})\in \R^{n_1\times \cdots \times n_d}$ and
$$
\vec{i}=\vec{j}\vec{k}
$$
be a partition of index $\vec{i}$, where
$$
\begin{array}{l}
\vec{j}=\vec{j}_1,\cdots,\vec{j}_r=Id(j_1,\cdots,Jr;n_{j_1},\cdots,n_{j_r}),\\
\vec{k}=\vec{k}_1,\cdots,\vec{k}_s=Id(k_1,\cdots,k_s;n_{k_1},\cdots,n_{k_s}),\\
\end{array}
$$
and $r+s=d$. Set $\a=\prod_{i\in {\bf j}}n_i$ and $\b=\prod_{i\in {\bf k}}n_i$,  a matrix
\begin{align}\label{2.3.2}
A=M^{\vec{j}\times \vec{k}}({\cal A})
\end{align}
is called a matrix expression of ${\cal A}$ if $A$ is obtained by arranging the elements of ${\cal A}$ in such a way that
the rows are in the order of $\vec{j}$ and the columns are in the order of $\vec{k}$.
\item[(ii)] Assume $\vec{j}=\vec{i}$ and $\vec{k}=\emptyset$. Then
\begin{align}\label{2.3.2}
V({\cal A}):=M^{\vec{i}\times \emptyset} ({\cal A})\in \R^{ {\bf n}}
\end{align}
is called the vector form of ${\cal A}$, where ${\bf n}=\prod_{i=1}^dn_i$.
\end{itemize}
\end{dfn}

\begin{exa}\label{e2.3.3}

Consider a hypermatrix ${\cal A}=(a_{i,j,k})\in \R^{2\times 3\times 5}$.
\begin{itemize}
\item[(i)]
$$
M^{\vec{i}\vec{j}\times \vec{k}}({\cal A})=
\begin{bmatrix}
a_{1,1,1}&a_{1,1,2}&a_{1,1,3}&a_{1,1,4}&a_{1,1,5}\\
a_{1,2,1}&a_{1,2,2}&a_{1,2,3}&a_{1,2,4}&a_{1,2,5}\\
a_{1,3,1}&a_{1,3,2}&a_{1,3,3}&a_{1,3,4}&a_{1,3,5}\\
a_{2,1,1}&a_{2,1,2}&a_{2,1,3}&a_{2,1,4}&a_{2,1,5}\\
a_{2,2,1}&a_{2,2,2}&a_{2,2,3}&a_{2,2,4}&a_{2,2,5}\\
a_{2,3,1}&a_{2,3,2}&a_{2,3,3}&a_{2,3,4}&a_{2,3,5}\\
\end{bmatrix}
$$
\item[(ii)]
$$
M^{\vec{k}\vec{i}\times \vec{j}}({\cal A})=
\begin{bmatrix}
a_{1,1,1}&a_{1,2,1}&a_{1,3,1}\\
a_{2,1,1}&a_{2,2,1}&a_{2,3,1}\\
a_{1,1,2}&a_{1,2,2}&a_{1,3,2}\\
a_{2,1,2}&a_{2,2,2}&a_{2,3,2}\\
a_{1,1,3}&a_{1,2,1}&a_{1,3,3}\\
a_{2,1,3}&a_{2,2,1}&a_{2,3,3}\\
a_{1,1,4}&a_{1,2,4}&a_{1,3,4}\\
a_{2,1,4}&a_{2,2,4}&a_{2,3,4}\\
a_{1,1,5}&a_{1,2,5}&a_{1,3,5}\\
a_{2,1,5}&a_{2,2,5}&a_{2,3,5}\\
\end{bmatrix}.
$$
\item[(iii)]
$$
\begin{array}{l}
V({\cal A})=[a_{1,1,1},a_{1,1,2},a_{1,1,3},a_{1,1,4},a_{1,1,5},\\
~~a_{1,2,1},a_{1,2,2},a_{1,2,3},a_{1,2,4},a_{1,2,5},\\
~~a_{1,3,1},a_{1,3,2},a_{1,3,3},a_{1,3,4},a_{1,3,5},\\
~~[a_{2,1,1},a_{2,1,2},a_{2,1,3},a_{2,1,4},a_{2,1,5},\\
~~a_{2,2,1},a_{2,2,2},a_{2,2,3},a_{2,2,4},a_{2,2,5},\\
~~a_{2,3,1},a_{2,3,2},a_{2,3,3},a_{2,3,4},a_{2,3,5}]^{\mathrm{T}}.\\
\end{array}
$$
\end{itemize}
\end{exa}

\section{Matrix Set vs Matrix Form}

 Consider a hypermatrix with order $2$, say, ${\cal A}=(a_{i,j})\in \R^{m\times n}$. It is worth to emphasize that
${\cal A}$  differs from a matrix. In fact it has $6$ different matrix expressions as
$$
\begin{array}{l}
A_1=M^{\{i,j\}\times \emptyset}({\cal A}) \in \R^{mn},\\
A_2=M^{\{j,i\}\times \emptyset}({\cal A}) \in \R^{mn},\\
A_3=M^{\{i\}\times \{j\}})({\cal A})\in {\cal M}_{m\times n},\\
A_4=M^{\emptyset \times \{i,j\}}({\cal A})=A_1^{\mathrm{T}} \in \R^{mn},\\
A_5=M^{\emptyset \times \{j,i\}}({\cal A})=A_2^{\mathrm{T}} \in \R^{mn},\\
A_6=M^{\{j\}\times \{i\}})({\cal A})=A_3^{\mathrm{T}}\in {\cal M}_{n\times m},\\
\end{array}
$$
When we consider the matrix properties, a special matrix expression must be pre-assigned!

Even when an order $1$ hypermatrix, it has two matrix forms as row vector and column vector. To distinct the elements in a matrix and the spacial matrix form of a set of data, we call a set of ordered data the {\bf matrix set}, and the special matrices constructed from this data set the {\bf matrix form}.

It is worth to emphasize that a hypermatrix is a ``matrix set" but not a ``matrix form". It is easy to verify that a hypermatrix of order $d$ has $(d+1)!$ different matrix forms. To some extent, the name  ``hypermatrix" is  confusing. It is a matrix set but not a matrix form. The objects discussed in Linear Algebra or Matrix Theory are matrices. Only when a special matrix expression of a hypermatrix is pre-assigned, the results from matrix theory can then be applied. Ask a hypermatrix to have matrix properties is ridiculous.

The following property is matrix form expression independent.

\begin{dfn}\label{d2.3.4}\cite{lim13} Given ${\cal A}=(a_{i_1,\cdots,i_d})\in \R^{n_1\times \cdots\times n_d}$ and a permutation $\sigma\in {\bf S}_d$.
The $\sigma$-transpose of ${\cal A}$, denoted by ${\cal A}^{\sigma}$ is defined by
\begin{align}\label{2.3.3}
{\cal A}^{\sigma}:=(a_{\sigma(1),\cdots,\sigma(d)}).
\end{align}
\end{dfn}

The vector form of a hypermatrix could be one of the most convenient matrix forms of  a hypermatrix.

A natural question is: how to get $V({\cal A}^{\sigma}$ from $V({\cal A})$? It will be seen in the sequel that this is a key issue in the further discussion.

Consider $\R^{n_1\times \cdots\times n_d}$ and   a permutation $\sigma\in {\bf S}_d$. Construct a matrix set of column vector as
$$
D_{\sigma}:=\left\{\d^{j_1}_{n_{\sigma(1)}}\d^{j_2}_{n_{\sigma(2)}}\cdots \d^{j_d}_{n_{\sigma(d)}}\;|\; j_k\in [1,n_{\sigma(k)}],\; k\in [1,d]\right\}.
$$

Then arrange the vectors in $D$ is the order as
$$
Id(j_{\sigma^{-1}(1)},j_{\sigma^{-1}(2)},\cdots,j_{\sigma^{-1}(d)};n_1,n_2,\cdots,n_d\}.
$$
The resulting matrix is denoted by
$$
W^{\sigma}_{[n_1\times n_2\times \cdots \times n_d]}
$$
is called the $sigma$-permutation matrix.

\begin{exa}\label{e2.3.5} Consider $\R^{3\times 4\times 2}$.
Assume $\sigma=(1,3,2)$. Then
$$
\begin{array}{l}
D_{\sigma}=\left\{\d^{j_1}_{n_{\sigma(1)}} \d^{j_2}_{n_{\sigma(2)}}\d^{j_3}_{n_{\sigma(3)}}\right\}\\
=\left\{\d^{j_1}_{n_{3}} \d^{j_2}_{n_{1}}\d^{j_3}_{n_{2}}\right\}
=\left\{\d^{j_1}_{2} \d^{j_2}_{3}\d^{j_3}_{4}\right\}\\
\end{array}
$$
Note that
$$
\sigma^{-1}(3)=1,~\sigma^{-1}(2)=3,~\sigma^{-1}(1)=2,
$$
we let $j_1$ run from $1$ to $2$ first, then $j_3$ run from $1$ to $4$, and finally $j_2$ run from $1$ to $3$. Then the permutation matrix is obtained as

$$
\begin{array}{l}
W^{\sigma}_{[3,4,2]}=[\d_2^1\d_3^1\d_4^1,\d_2^2\d_3^1\d_4^1,
\d_2^1\d_3^1\d_4^2,\d_2^2\d_3^1\d_4^2,\d_2^1\d_3^1\d_4^3,\\
~~\d_2^2\d_3^1\d_4^3,\d_2^1\d_3^2\d_4^1,\d_2^2\d_3^2\d_4^1,\d_2^1\d_3^2\d_4^2,\d_2^2\d_3^2\d_4^2,
\d_2^1\d_3^2\d_4^3,\\
~~\d_2^2\d_3^2\d_4^3,\d_2^1\d_3^3\d_4^1,\d_2^2\d_3^3\d_4^1,
\d_2^1\d_3^3\d_4^2,\d_2^2\d_3^3\d_4^2,\d_2^1\d_3^3\d_4^3,\d_2^2\d_3^3\d_4^3].\\
=\d_{24}[1,13,2,14,3,15,4,16,5,17,6,18,\\
~~7,19,8,20,9,21,10,22,11,23,12,24].
\end{array}
$$
\end{exa}

Using permutation matrix, the aforementioned question can be answered.

\begin{prp}\label{p2.3.6} \cite{che12} The $V({\cal A}^{\sigma})$ can be calculated by the following formula.
\begin{align}\label{2.3.4}
V({\cal A}^{\sigma})=W^{\sigma}_{[n_1\times\cdots\times n_d]}V({\cal A}).
\end{align}
\end{prp}

\begin{rem}\label{r2.3.7} Formula (\ref{2.3.4}) is a fundamental tool in our later study. In fact, it says that if we has a matrix set
as
$$
D=\{a_{i_1,\cdots.i_d}\;|\;i_k\in [1,n_k],~k\in [1,d]\},
$$
We can arrange the elements into a vector using different orders. The order $Id(i_1,\cdots,i_d;n_1,\cdots,n_d)$ is called the natural order.
Using the natural order we have a vector $x$. Then we use another order as $Id(i_{\sigma(1)},\cdots,i_{\sigma(d)}; n_{\sigma(1)},\cdots,n_{\sigma(d)})$ to get  another vector $y$. Then we have
$$
y=W^{\sigma}_{[n_1\times\cdots\times n_d]}x.
$$
\end{rem}

\section{KPD of Hypervectors}

\subsection{Monic Decomposition Algorithm}

Monic Decomposition Algorithm (MDA) was firstly proposed by us in \cite{chepr}. It consists of a set of linear mapping from higher dimensional Euclidian space to its factor-dimensional subspaces. We firts give a detailed review on MDA. And then use it for the decomposition of vectors.
%

\begin{prp}\label{pkd.1.2.2}
Assume $x=\ltimes_{i=1}^dx_i$, where $x_i\in \R^{n_i}$, $i\in [1,d]$, and $e(x)=e$, and $e(x_i)=e_i$, $i\in [1,d]$.
Set $j:=e$ and $(i_1,\cdots,i_d):=(e_1,e_2,\cdots,e_d)$. Then $e$ and $(e_1,\cdots,e_d)$ satisfy (\ref{kd.1.1.7})-(\ref{kd.1.1.8}).
\end{prp}

\noindent{\it Proof.} Using $\vec{j}=\vec{i}_1\vec{i}_2\cdots \vec{i}_d$ to label the elements of $x$. It is easy to verify the fact that $h_0(x)=\prod_{i=1}^d h_0(x_i)$. That is,
\begin{align}\label{kd.1.2.1}
x_{e}=(x_1)_{e_1}(x_2)_{e_2}\times (x_d)_{e_d}.
\end{align}
Hence, the index $j=d$ corresponds the product of indexes $\vec(i)_k=e_k$, $k\in [i,d]$. Then the conclusion follows immediately.
\hfill $\Box$

\begin{dfn}\label{dkd.1.2.3} (KPD for vectors) Consider $x\in \R^{{\bf n}}$, where ${\bf n}=\prod_{i=1}^dn_i$.
$x$ is said to be decomposable w.r.t. $n_1\times \cdots\times n_d$, if there exits
$x_i\in \R^{n_i}$, $i\in [1,d]$ such that
\begin{align}\label{kd.1.2.2}
x=\ltimes_{i=1}^dx_i.
\end{align}
\end{dfn}

It is easy to verify the following result \cite{che19}, which shows if there exists a decomposition, then it is unique up to a set of constant product coefficients.

\begin{prp}\label{pkd.1.2.4} Assume $0\neq x=\ltimes_{i=1}^d x_i=\ltimes_{i=1}^dz_i$, then
\begin{align}\label{kd.1.2.3}
z_i=\mu_i x_i,\quad i\in [1,d],
\end{align}
where $\prod_{i=1}^d \mu_i=1$.
\end{prp}

Using Proposition \ref{pkd.1.2.2}, for a given vector $x$ with its known head index $e=e(x)$, if it is decomposable then all $e_i=e(x_i)$ are well known. Then we can define a set of decomposition mappings as
\begin{align}\label{kd.1.2.4}
\begin{array}{l}
\Xi^e_{(i;d)}:=[\d_{n_1}^{e_1}]^{\mathrm{T}} \otimes \cdots  \otimes [\d_{n_{i-1}}^{e_{i-1}}]^{\mathrm{T}}\otimes I_{n_i}\otimes\\
~~[\d_{n_{i+1}}^{e_{i}+1}]^{\mathrm{T}} \otimes \cdots  \otimes [\d_{n_{r}}^{e_{r}}]^{\mathrm{T}},\quad i\in [1,d].
\end{array}
\end{align}

The following Proposition shows how to calculate the component vectors from a hypervector. It is called the monic decomposition algorithm (MDA). 

\begin{prp}\label{pkd.1.2.5}
Consider $x\in \R^{{\bf n}}$, where ${\bf n}=\prod_{i=1}^dn_i$, and $e(x)=e$, $h_0(x)=a$.
Assume $x$ is decomposable w.r.t., $n_1\times \cdots\times n_d$, and set $x_0=\frac{1}{a}x$. Define
\begin{align}\label{kd.1.2.5}
x_i:=\Xi^e_{(i;d)}(x_0),\quad i\in [1,d],
\end{align}
then all $x_i$, $i\in [1,d]$ are monic. Moreover,
\begin{align}\label{kd.1.2.6}
x=a\ltimes_{i=1}^dx_i.
\end{align}
\end{prp}

\noindent {\it Proof} Note that for tow column vectors $\xi$ and $\eta$ $\xi \ltimes\eta=\xi\otimes \eta$. Then we have
$$
\begin{array}{l}
[(\d_{n_k}^{e_k}]^{\mathrm{T}}x_k=1,\quad k\in [1,d]~\mbox{and}~k\neq i ,\\
I_i\otimes x_i=x_i.
\end{array}
$$
Using the formula that $(A\otimes B)(C\otimes D)=(AC)\otimes (BD)$, Equations (\ref{kd.1.2.5}) follows.
\hfill $\Box$

\begin{rem}\label{rkd.1.2.501}
\begin{itemize}
\item[(i)] If $x$ is not decomposable, then we can still use the decomposition mappings to get the projections $x_k$, $k\in [1,d]$. Now of course, (\ref{kd.1.2.6}) does not hold. Hence, MDA not only provides an algorithm to calculate the component vectors of a hypervector, but also provides a way to verify whether a vector is decomposable.

\item[(ii)] Even when $x$ is not decomposable, the $x_k$, $k\in [1,d]$ can be considered as  component vectors of an approximated decomposition of $x$. Because the decomposition mappings are linear, when $x$ is close to a decomposable hypervector, the
$\ltimes_{k=1}^sx_k$ can be thought as an approximated decomposition.
\end{itemize}
\end{rem}

\begin{exa}\label{ekd.1.2.6} Given a hypervector $x\in \R^{3\ltimes 4\ltimes 2}$, with
$$
x_{14}=2,~x_{16}=1,~x_{22}=-2,~x_{24}=-1,
$$
and the other entries are $0$.
We use two methods to find the components of $x$.
\begin{itemize}
\item[(i)] Simultaneous Decomposition:
Set
$$
x=x_1\ltimes x_2\ltimes x_3.
$$
Since $e(x)=14$, and $h_0(x)=2$, then
$x_0=x/2$, and
it is easy to calculate that
$$
\begin{array}{l}
e_1:=e(x_1)=13-[13/2]*2+1=2,\\
e_2:=e(x_2)=6-[6/4]*4+1=3,\\
e_3=1+1=2.
\end{array}
$$
Then we construct the decomposition mappings as
$$
\begin{array}{l}
\Xi^e_{(1;3)}=I_3\otimes(\d_8^6)^{\mathrm{T}}:=(\xi_{i,j})\in {\cal M}_{3\times 24},\\
~~\mbox{where}~\xi_{1,6}=1,~\xi_{2,14}=1,~\xi_{3,22}=1,~\mbox{all otheres are}~0.\\
\Xi^e_{(2;3)}=(\d_{3}^2)^{\mathrm{T}}\otimes I_4\otimes(\d_2^2)^{\mathrm{T}}:=(\xi_{i,j})\in {\cal M}_{4\times 24},\\
 ~~\mbox{where}~\xi_{1,10}=1,~\xi_{2,12}=1,~\xi_{3,14}=1,~\xi_{4,16}=1,\\
 ~~~\mbox{all otheres are}~0.\\
\Xi^e_{(3;3)}=(\d_{12}^7)^{\mathrm{T}} \otimes  I_2:=(\xi_{i,j})\in {\cal M}_{2\times 24},\\
 ~~\mbox{where}~\xi_{1,13}=1,~\xi_{2,14}=1,~\mbox{all otheres are}~0.\\
\end{array}
$$
It follows that
$$
\begin{array}{l}
x_1=\Xi^2_{(1;3)}x_0=(0,1,-1)^{\mathrm{T}},\\
x_2=\Xi^2_{(2;3)}x_0=(0,0,1,0.5)^{\mathrm{T}},\\
x_3=\Xi^2_{(3;3)}x_0=(0,1)^{\mathrm{T}},\\
\end{array}
$$
Finally, we have
$$
x=h_0(x)x_1\ltimes x_2\ltimes x_3=2x_1\ltimes x_2\ltimes x_3.
$$
\item[(ii)] Step by Step  Decomposition:

 Step 1: Consider $x\in \R^{3\ltimes 8}$ as in (i). Decompose
 $x=h_0x_1\ltimes z$.
Similarly to (i), we can get
 $$
 \begin{array}{l}
 x_1=(0,1,-1)^{\mathrm{T}},\\
 z=(0,0,0,0,0,1,0,0.5)^{\mathrm{T}}.
\end{array}
$$
Next, decompose $z$ yields
$$
\begin{array}{l}
x_2=\Xi^2_{(2;3)}x_0=(0,0,1,0.5)^{\mathrm{T}},\\
x_3=\Xi^2_{(3;3)}x_0=(0,1)^{\mathrm{T}},\\
\end{array}
$$
\end{itemize}
\end{exa}

\subsection{Approximate Decomposition}

Given an $x\in \R^{{\bf n}}$ and assume ${\bf n}=\prod_{i=1}^dn_i$. Then its head index $e(x)$ and head value $h_0$ are known. Setting $x_0=x/h_0$ and using decomposition mapping we can always obtained
$x_i=\Xi^e_{(i;d)}x_0,\quad i\in [1,d]$. Then MDA can be used to verify  whether $x$ is decomposable.

\begin{exa}\label{ekd.1.3.2} Consider Example \ref{ekd.1.2.6} again. Let $x\in \R^{3\ltimes 4\ltimes 2}$ be the same as in Example \ref{ekd.1.2.6} except $x_{16}=4$.
Then it is a straightforward computation to obtain that
\begin{align}\label{kd.1.3.2}
\begin{array}{l}
x_1=(0,1,-1)^{\mathrm{T}},\\
x_2=(0,0,1,2)^{\mathrm{T}},\\
x_3=(0,1)^{\mathrm{T}}.
\end{array}
\end{align}
Define the square error as
\begin{align}\label{kd.1.3.3}
E=(x_0-x_1x_2x_3)^{\mathrm{T}}(x_0-x_1x_2x_3).
\end{align}
Then it is easy to calculate that for the solution in (\ref{kd.1.3.2}) the square  error is $E=2.25$. Hence the $x$ is not decomposable.
\end{exa}

As aforementioned that when $x$ is not decomposable, the set of components obtained by decomposition mappings form an approximate decomposition of $x$. It might be expected that this set of components is the least square error approximation.  Unfortunately, a numerical computation shows that this guess is wrong.

To get the least square approximated KPD for non-decomposable vectors, we provide the following algorithm. 

\begin{alg}\label{akd.1.3.3} Given $x\in \R^{{\bf n}}$, with ${\bf n}=\prod_{i=1}^dn_i$. The least square approximated KPD
can be obtained as follows.
\begin{itemize}
\item[(i)] Calculate $e=e(x)$, $h_0=x_e$, $x_0=x/h_0$, and then using Proposition \ref{pkd.1.2.2} to calculate $e_i$, $i\in [1,d]$.

\item[(ii)] Define candidate component vectors as
\begin{align}\label{kd.1.3.4}
x_i=(\underbrace{0,\cdots,0}_{e_i-1}, 1,u^i_1,\cdots,u^i_{n_i-e_i})^{\mathrm{T}},\quad i\in [1,d].
\end{align}

\item[(iii)] Using MDA to get the initial values of the iteration as
\begin{align}\label{kd.1.3.5}
x_i(0)=\Xi^e_{(i;d)}x_0,\quad i\in [1,d].
\end{align}
(\ref{kd.1.3.5}) is equivalent to initialize $u=u(0)$, where $u=\{u^i_j\;|\; j\in [1,n_i-e_i],~i\in [1,d]\}$.

\item[(iv)] Consider the square error
\begin{align}\label{kd.1.3.6}
E(u)=(x_0-\ltimes_{i=1}^dx_i)^{\mathrm{T}}(x_0-\ltimes_{i=1}^dx_i).
\end{align}
Its gradients are
\begin{align}\label{kd.1.3.7}
D_{i,j}:=\frac{\pa E(u)}{\pa u^i_j},\quad j\in [1,n_i-e_i],~i\in [1,d].
\end{align}

Using gradient descent method, we choose a step length $h>0$ and define the iteration as
\begin{itemize}
\item[] Step 1:
Set $u^i_j(0)$ as in (\ref{kd.1.3.5}).
\item[] Step 2:
Calculate $E(0)$ by (\ref{kd.1.3.6}).
\item[] Step k:
Set
$$
u^i_j(k+1)=u^i_j(k)-D_{i,j}*h, \quad j\in [1,n_i-e_i],~i\in [1,d].
$$
Calculate
$$
E(k+1)=E(u(k+1)).
$$
If $E(k+1)<E(k)$, set $u^i_j(0)=u^i_j(k+1)$ and go to Step 2.

Else Stop.
\end{itemize}
\end{itemize}
\end{alg}

\begin{prp}\label{pkd.1.3.4} The algorithm \ref{akd.1.3.3}  converges to the least square error solution.
\end{prp}
\noindent{\it Proof.} It is easy to see that for a fixed $u^s_t$, assume $u^i_j$, $(i,j)\neq (s,t)$ are fixed, then the $E(u)$ is a quadratic form of $u_{i,j}$. Hence $E(u)$ is a higher dimensional paraboloid. The minimum value point is unique, and it is the only  stationary point. Hence  the gradient descent algorithm must converge to the minimum value.
\hfill $\Box$

\begin{exa}\label{ekd.1.3.5}  Consider Example \ref{ekd.1.3.2} again.
Using Algorithm \ref{akd.1.3.3},  we can get  the candidate least square approximated KPD as
$$
\begin{array}{l}
x_1=(0,1,u),\\
x_2=(0,0,1.v),\\
x_3=(0,1).
\end{array}
$$
The initial value is $(u,v)=(-1,2)$.
The square error is:
$$
E=(v-2)^2+(u+1)^2+(uv+0.5)^2.
$$
It follows that
$$
\begin{array}{l}
D1:=\frac{\pa E}{\pa u}=2(u+1+v(uv+0.5)),\\
D2:=\frac{\pa E}{\pa v}=2(v-2+u(uv+0.5)).\\
\end{array}
$$
Choosing $h=0.05$, initial value  $(u,v)=(-1,2)$, and using the gradient descent algorithm,  after 100 iterations we have
$$
u^*=-0.4306,\quad v^*=1.8688,
$$
and the square error is
$$
E^*=0.4343,
$$
comparing with $E(0)=2.25$.
\end{exa}

\subsection{Finite Sum Decomposition}

Next, we consider the finite sum KDP, that is, find small number of vectors  $(x_1^k,x_2^k)$ such that
\begin{align}\label{kd.1.4.1}
x=\dsum_{k=1}^r x_1^k\ltimes x_2^k.
 \end{align} 

\begin{rem}\label{rkd.1.4.0} For statement ease, this subsection considers only the two factor decomposition. In fact, when the KPD of hypervectors is considered, there is no difference between two-factor KPD and $n$-factor KPD. This fact can be seen in Examples \ref{ekd.1.2.6} and \ref{ekd.1.3.2}. But when the technique is extended to  KPD of matrices and KPD of hypermatrices, the algorithm will be developed for two-factor KPD only. Hence for  matrix and hypermatrix cases, the $n$ factor KPD must be executed by two-factor KPDs step by step.
\end{rem}

We give the following algorithm to realize  multi-term KDP.

\begin{alg}\label{akd.1.4.1}  Given $x\in \R^{{\bf n}}$, with ${\bf n}=\prod_{i=1}^dn_i$, and a pre-assigned error boundary $\epsilon >0$.
\begin{itemize}
\item[] Step 1:
Set $x(0)=x$, using Algorithm \ref{akd.1.3.3} to get the least square approximated KDP
\begin{align}\label{kd.1.4.2}
x(0)\approx a_1\ltimes_{i=1}^dx^1_i,
\end{align}
and
\begin{align}\label{kd.1.4.3}
x(1)=x(0)- a_1\ltimes_{i=1}^dx^1_i.
\end{align}
\item[] Step k: (Iteration)
Using Algorithm \ref{akd.1.3.3} to get least square KDP
\begin{align}\label{kd.1.4.4}
x(k-1)\approx a_k\ltimes_{i=1}^dx^k_i,
\end{align}
and
\begin{align}\label{kd.1.4.5}
x(k)=x(k-1) - a_k\ltimes_{i=1}^dx^k_i.
\end{align}
When $\|x(k)\|=0$ , step.
\end{itemize}
\end{alg}

\begin{prp}\label{pkd.1.4.101} The algorithm  \ref{akd.1.4.1} will stop at finite times.
\end{prp}

\noindent{\it Proof.} Assume the head index $e(x)=e$ and head value is $h_0(x)=a$. Using MDA the approximated decomposition is
$\tilde{x}:=a\ltimes_{i=1}^dx_i$. It is easy to see that $e(\tilde{x})=e$ and $h_0(\tilde{x})=a$. It follows that
$$
e(x-\tilde{x})>e(x).
$$
Hence after finite steps either the $x-\tilde{x}=0$  or $e(x-\tilde{x})=n$. It is obvious that  if $e(x)=n$, $x$ is decomposable.
\hfill $\Box$

\begin{exa}\label{ekd.1.4.2}  Consider Example \ref{ekd.1.3.5} again.
From Example \ref{ekd.1.3.5} we have
\begin{align}\label{kd.1.4.5}
x_0:=x\approx h^1_0x^1_1x^1_2x^1_3,
\end{align}
where $h^1_0=2$, and
$$
\begin{array}{l}
x^1_1=[0,1,-0.4306]^{\mathrm{T}},\\
x^1_2=[0,0,1,1.8688]^{\mathrm{T}},\\
x^1_3=[0,1]^{\mathrm{T}}.\\
\end{array}
$$
Set
$$
x_1=x_0-h^1_0x^1_1x^1_2x^1_3,
$$
and using the least square approximating KPD to $x_1$ yields
$h^2_0=0.2624$, and
$$
\begin{array}{l}
x^2_1=[0,1,2.3221]^{\mathrm{T}},\\
x^2_2=[0,0,0,1]^{\mathrm{T}},\\
x^2_3=[0,1]^{\mathrm{T}}.\\
\end{array}
$$
Set
$$
x_2=x_1-h_1x^2_1x^2_2x^2_3.
$$
It is easy to verify that $x_2$ is decomposable, which can be decomposed as
$h^3_0=-1.1388$, and
$$
\begin{array}{l}
x^3_1=[0,1,1]^{\mathrm{T}},\\
x^3_2=[0,0,1,0]^{\mathrm{T}},\\
x^2_3=[0,1]^{\mathrm{T}}.\\
\end{array}
$$
We conclude that
$$
x=h^1_0 x^1_1x^1_2x^1_3+h^2_0 x^2_1x^2_2x^2_3+ h^3_0x^3_1x^3_2x^3_3.
$$
\end{exa}

\begin{rem}\label{rkd.1.4.3}
\begin{itemize}
\item[(i)] It was proved in \cite{van93} that a matrix always has explicit finite sum KPD. In vector case it can be seen that  the maximum time is  $n=\dim(x)$.
Since in our algorithm each time the vector index $e$ is increasing at least by $1$, so at most  after $n$ times the zero error finite sum KPD is obtained.

\item[(ii)] The Algorithm \ref{akd.1.4.1} can stop at any time when the $\|x(k)\|<\epsilon$, which is to your satisfaction.

\item[(iii)] Since
$$x=\ltimes_{i=1}^dx_i=x_1\ltimes (x_2\ltimes (\cdots \ltimes x_d)\cdot),
$$
it is obvious that the KPD can be done via $2$-component decompositions step by step as
\begin{align}\label{kd.1.4.6}
\begin{array}{l}
x=x_1\ltimes x^2,\\
 x^2=x_2\ltimes x^3,\\
 \vdots,\\
 x^{k-1}=x_{k-1}\ltimes x_k.
 \end{array}
\end{align}

 But when approximate KPD is considered, can the least square approximate KPD  be obtained via   two component decompositions step by step as in (\ref{kd.1.4.6}) (with ``$=$"being replaced by ``$\approx$")?
Based on numerical computations our conjecture is ``Yes". But if the conjecture does not hold, the step by step approximated KPD might be inferior to the simultaneous approximation.
\end{itemize}
\end{rem}

\section{KPD of Matrices}

\begin{dfn}\label{dkd.2.1} Let $N\in {\cal M}_{{\bf m}\times {\bf n}}$ with ${\bf m}=\prod_{i=1}^dm_i$, and ${\bf n}=\prod_{i=1}^dn_i$.  The KPD problem is:  finding $B_i\in {\cal M}_{m_i\times n_i}$, $i\in [1,d]$, such that
\begin{align}\label{kd.2.1}
N=B_1\otimes \cdots \otimes B_d.
\end{align}
\end{dfn}\index{KPD}

Here we consider only the case when $d=2$. Given $N=(n_{u,v})\in {\cal M}_{mp\times nq}$, our goal is to find $A=(a_{i,j})\in {\cal M}_{m\times n}$ and $B=(b_{r,s})\in {\cal M}_{p\times q}$, such that
$$
N=A\otimes B.
$$

First, we calculate $V_c(A\otimes B)$: Note that
$$
A\otimes B=[a^{1}\otimes b^{1},\cdots, a^{1}\otimes b^{q}, a^{2}\otimes b^{1},\cdots,a^{n}\otimes b^{q}],
$$
where $a^{i}=\Col_i(A)$, $i\in [1,n]$, $b^{i}=\Col_i(B)$, $i\in [1,q]$.
It follows that
$$
V_c(A\otimes B)=\begin{bmatrix}
a^{1}\otimes b^{1}\\
a^{1}\otimes b^{2}\\
\cdots\\
a^{1}\otimes b^{q}\\
a^{2}\otimes b^{1}\\
\vdots\\
a^{n}\otimes b^{q}\\
\end{bmatrix}.
$$
Setting
$$
c^{i,j}_{r,s}:=a_{i,j}b_{r,s},
$$
then a straightforward computation shows that the elements of $V_c(A \otimes B)$ are
$$
D=\{c^{i,j}_{p,q}\;|\;i\in [1,m], j\in [1,n],r\in[1,p],s\in [1,q]\},
$$
and they are arranged in the order of
$$
Id(j,s,i,r;n\times q\times m\times p)=\vec{j}\vec{s}\vec{i}\vec{r}.
$$

On the other hand $V(A)=V_c(A)$ with elements in
$$
D_A=\{a_{i,j}\;|\;i\in [1,m],j\in [1,n]\},
$$
and they are arranged in the order of $Id(j,i;n,m)=\vec{j}\vec{i}$.

$V(B)=V_c(B)$ with elements in
$$
D_B=\{b_{r,s}\;|\;r\in [1,p],s\in [1,q]\},
$$
and they are arranged in the order of $Id(s,r;q,p)=\vec{s}\vec{r}$.
Then $V(A)\ltimes V(B)$ has elements in $D$ with the order as
$$
Id(j,i,s,r;n,m,q,p)=\vec{j}\vec{i}\vec{s}\vec{r}.
$$
Now to convert $\vec{j}\vec{i}\vec{s}\vec{r}$ to $\vec{j}\vec{s}\vec{i}\vec{r}$, a swap matrix $W_{[m,q]}$ can be used to construct
\begin{align}\label{kd.2.101}
\Psi:=I_n\otimes W_{[m,q]}\otimes I_p.
\end{align}
Then it is clear that
\begin{align}\label{kd.2.2}
V(A\otimes B)=\Psi V(A)\ltimes V(B).
\end{align}
Since $\Psi$ is an orthogonal matrix, (\ref{kd.2.3}) is equivalent to
\begin{align}\label{kd.2.3}
\Psi^{\mathrm{T}}V(A\otimes B)=V(A)\ltimes V(B).
\end{align}

This leads to the following result immediately.

\begin{thm} \label{tkd.2.2} Given $N\in {\cal M}_{mp\times nq}$. $N$ is paired Kronecker product decomposable w.r.t. $(m\times n,~p\times q)$, if and only if, $\Psi^{\mathrm{T}}V(N)$ is decomposabe w.r.t. $mn\times pq$, where $\Psi$ is defined in (\ref{kd.2.101}).
\end{thm}

We use an example to descripe the decomposition process.

\begin{exa}\label{ekd.2.3}
Given
$$
A=\begin{bmatrix}
0&0&0&1&2&-1\\
0&0&0&1&0&-2\\
-1&-2&1&1&2&-1\\
-1&0&2&1&0&-2\\
\end{bmatrix}\in {\cal M}_{4\times 6}.
$$
Consider its paired KPD w.r.t. $2\times 2,2\times 3$.

Using (\ref{kd.2.101}), we construct
$$
\begin{array}{l}
\Psi=I_2\otimes W_{[2,3]}\otimes I_2=
\d_{24}[1,2,7,8,3,4,9,10,\\
~~5,6,11,12,13,14,19,20,15,16,21,22,17,18,23,24].
\end{array}
$$
Then
$$
\begin{array}{l}
V:=\Psi^{\mathrm{T}} V(A)=[0,0,0,0,0,0,-1,-1,\\
~~-2,0,1,2,1,1,2,0,-1,-2,1,1,2,0,-1,-2]^{\mathrm{T}}.
\end{array}
$$
Hence $e(V)=7$ and $h_0(V)=-1$. Set $x=V/h_0$ and then use MDA on $x$.
First it is easy to calculate the $e_1=2$ and $e_2=1$. Then the decomposition mappings can be constructed as
$$
\begin{array}{l}
\Xi^e_{(1;2)}=I_4\otimes [\d_6^1]^{\mathrm{T}},\\
\Xi^e_{(2;2)}=[\d_4^2]^{\mathrm{T}}\otimes I_6.
\end{array}
$$
Hence we have
$$
\begin{array}{l}
x_1=\Xi^e_{(1;2)}x=(0,1,-1,-1)^{\mathrm{T}},\\
x_2=\Xi^e_{(2;2)}x=(1,1,2,0,-1,-2)^{\mathrm{T}}.\\
\end{array}
$$
It follows immediately that
$$
B=\begin{bmatrix}
0&-1\\
1&-1\\
\end{bmatrix}
$$
$$
C=\begin{bmatrix}
1&2&-1\\
1&0&-2\\
\end{bmatrix}
$$
It is easy to verify that
$$
A=h_0 B\otimes C=-B\otimes C.
$$
\end{exa}

\begin{rem}\label{rkd.2.3}
\begin{itemize}
\item[(i)] Consider the multi-fold decomposition $A=B_1\otimes\cdots\otimes B_d$.  Assume $A$ is decomposable, the it is obvious that it can be done via two-fold decomposition step by step.
\item[(ii)] When the approximate multi-fold decomposition is considered, of course, the step by step two-fold decomposition is still applicable. But if the conjecture in Remark \ref{rkd.1.4.3} is not true, the solution may be inferior to the least square error solution.
\item[(iii)] Unlike in vector case, in matrix case the simultaneous multi-fold decomposition seems difficult.
\end{itemize}
\end{rem}

Next, we consider the approximate KPD. It is immediate consequence of Theorem \ref{tkd.2.2} that if the exact decomposition does not exist, then, similarly to vector case, we can use the components obtained from decomposition mappings as the initial value and use the gradient descent method to find the least square error approximate decomposition. We use an example to demonstrate this.


\begin{exa}\label{ekd.2.7}
Consider
$$
A=\begin{bmatrix}
0&0&0&1&2&-1\\
0&0&0&1&0&-2\\
2&-2&1&1&2&-1\\
-1&0&2&1&0&-2\\
\end{bmatrix}\in {\cal M}_{4\times 6}.
$$
It is similar to the one in Example \ref{ekd.2.3} only the $A(3,1)$ is changed from $-1$ to $2$.
We have $e(V(A))=7$ and $h_0(V(A))=2$.
Let $x_0=\frac{1}{h_0}V(A)$. Then we can calculate
$$
\begin{array}{l}
x:=\Psi^{\mathrm{T}} x_0=[0,0,0,0,0,0,1,-0.5,-1,0,0.5,1,\\
~~0.5,0.5,1,0,-0.5,-1,0.5,0.5,1,0,-0.5,-1].
\end{array}
$$
Similarly to Example   \ref{ekd.2.3}, we can have the decomposed components as
$$
\begin{array}{l}
x_1=(0,1,-1,-1)^{\mathrm{T}} ,\\
x_2=(1,-0.5,-1,0, 0.5, 1)^{\mathrm{T}} ,\\
\end{array}
$$
The square error is
$$
E=(x-x_1x_2)^{\mathrm{T}}(x-x_1x_2),
$$
which is
$$
E_0=11.25.
$$

Next, we assume the least square approximation is
$$
\begin{array}{l}
x_1=(0,1,u,v)^{\mathrm{T}} ,\\
x_2=(1,y_1,y_2,y_3, y_4, y_5)^{\mathrm{T}} ,\\
\end{array}
$$

The square error is calculated as
$$
\begin{array}{l}
E=(y_1+0.5)^2+(y_2+1)^2+(y_3)^2+(y_4-0.5)^2+(y_5-1)^2+\\
~~(u-0.5)^2+(uy_1-0.5)^2+(uy_2-1)^2+(uy_3)^2\\
~~+(uy_4+0.5)^2+(uy_5+1)^2+(v-0.5)^2+(vy_1-0.5)^2\\
~~+(vy_2-1)^2+(vy_3)^2+(vy_4+0.5)^2+(vy_5+1)^2.
\end{array}
$$
The gradients are
$$
\begin{array}{l}
\frac{\pa E}{\pa u}=2[(u-0.5)+y_1(uy_1-0.5)+y_2(uy_2-1)+y_3(uy_3)\\
~~+y_4(uy_4+0.5)+y_5(uy_5+1)],\\
\frac{\pa E}{\pa v}=2[(v-0.5)+y_1(vy_1-0.5)+y_2(vy_2-1)+y_3(vy_3)\\
~~+y_4(vy_4+0.5)+y_5(vy_5+1)],\\
\frac{\pa E}{\pa y_1}=2[(y_1+0.5)+u(uy_1-0.5)+v(vy_1-0.5)],\\
\frac{\pa E}{\pa y_2}=2[(y_2+1)+u(uy_2-1)+v(vy_2-1)],\\
\frac{\pa E}{\pa y_3}=2[(y_3)+u(uy_3)+v(vy_3)],\\
\frac{\pa E}{\pa y_4}=2[(y_4-0.5)+u(uy_4+0.5)+v(vy_4+0.5)],\\
\frac{\pa E}{\pa y_5}=2[(y_5-1)+u(uy_5+1)+v(vy_5+1)].\\
\end{array}
$$
Choosing $h=0.1$, setting
initial value
$$
\begin{array}{l}
u(0)=-1,~v(0)=-1; ~y_1(0)=-0.5,~y_2(0)=-1,\\
~y_3(0)=0,~y_4(0)=0.5,~y_5(0)=1,
\end{array}
$$
 and using the gradient descent algorithm,  after 248 iterations we have
$$
\begin{array}{l}
u^*=-0.5223,~v^*=-0.5223; ~y_1^*=-0.6614,~y_2^*=-1.3229,\\
~y_3^*=0,~y_4^*=0.6614,~y_5^*=1.3229.
\end{array}
$$
The corresponding square error is $E^*=2.8284$, comparing with $E_0= 11.25$.

We conclude that the least square error approximation is
$$
A\approx h_0 (B\otimes C)
$$
where $h_0=2$, and
$$
B=\begin{bmatrix}
0&-0.5223\\
1&-0.5223\\
\end{bmatrix}
$$
$$
C=\begin{bmatrix}
1&-1.3229&0.6614\\
-0.6614&0&1.3229\\
\end{bmatrix}
$$
 \end{exa}

\begin{rem}\label{rkd.2.8} Observing Example \ref{ekd.2.7} again. Of cause we can set $A_1:=A-h_0 (B\otimes C)$ and find the least square error approximation for $A_1$. Repeating this process finite times the finite sum KPD can be obtained. The procedure is exactly the same as we did in Example  \ref{ekd.2.7}. We skip the details here.
\end{rem}

\section{Kronecker Products of HyperMatrices}

\subsection{Form-based Kronecker Product}

 Consider ${\cal A}\in \R^{m\times n}$ and   ${\cal B}\in \R^{u\times v}$. Can we define a Kronecker product for them? In fact, the answer is ``No". Then the Kronecker Product of $A=(a_{i,j})\in {\cal M}_{m\times n}$ and $B=(b_{r,s})\in {\cal M}_{u\times v}$ is properly defined. Where the $A$ and $B$  are considered as two normal matrix expressions of ${\cal A}$ and ${\cal B}$ respectively.
 If we change the matrix expressions  $V({\cal A})$ and $V({\cal B})$, the corresponding Kronecker product is completely different from the previous one.

 From this perspective one sees easily that the Kronecker product of two order $2$ hypermatrices, as two matrix sets, can not be uniquely defined. While the Kronecker product of two matrices, as two matrix forms, is properly defined.

 Because of this, we may preassign two special matrix forms, say $A$ and $B$  for the hypermatrices ${\cal A}$ and ${\cal B}$ respectively. Then we define a form-based Kronecker product of ${\cal A}$ and ${\cal B}$ as
 \begin{align}\label{kd.3.1.0}
 {\cal A}\otimes {\cal B}:=A\otimes B.
 \end{align}
 In this way, many different kinds of Kronecker products can be defined.

 In fact, whole matrix theory is based on the fixed matrix forms of order $1$ and order $2$ hypermatrices,
 which are vectors and matrices respectively. That is, let ${\cal A}=(a_i)\in \R^{n}$. Then we use
 $x=M^{i\times \emptyset}({\cal A})$ as its default matrix form, which is a vector. Similarly,  let ${\cal A}=(a_{i,j}\in \R^{m\times n}$. Then we use
 $A=M^{i\times j}({\cal A})$ as its default matrix form, which is a matrix. Then the whole matrix theory is based on these default matrix forms. We may call these matrix forms as normal matrix form. When other forms are required, the $x^{\mathrm{T}}$, or $A^{ \mathrm{T}}$, $V_r(A)$, $V_c(A)$ etc. must be defined.

 To develop the matrix theory for hypermatrices, the (default) normal matrix form of order $d$ hypermatrices is very reasonable.
 In a normal form the operations and the numerical computations become much easier.

 For instance, consider a cubic matrix  ${\cal A}=(a_{k,i,j})\in \R^{m\times n\times s}$.  To investigate the t-product of cubic matrices and its applications, a normal form of cubic matrices is proposed as follows\cite{arz18,chen24}:
\begin{align}\label{kd.3.1.1}
A=\begin{bmatrix}
a^1_{1,1}&a^1_{1,2}&\cdots&a^1_{1,n}\\
a^1_{2,1}&a^1_{2,2}&\cdots&a^1_{2,n}\\
\cdots&~&~&~\\
a^1_{m,1}&a^1_{m,2}&\cdots&a^1_{m,n}\\
\vdots&~&~&~\\
a^s_{1,1}&a^s_{1,2}&\cdots&a^s_{1,n}\\
a^s_{2,1}&a^s_{2,2}&\cdots&a^s_{2,n}\\
\cdots&~&~&~\\
a^s_{m,1}&a^s_{m,2}&\cdots&a^s_{m,n}\\
\end{bmatrix}
=\begin{bmatrix}
A^{(1)}\\
A^{(2)}\\
\vdots\\
A^{(s)}
\end{bmatrix},
\end{align}
where $a^k_{i,j}=a_{k,i,j}$.

 \subsection{Outer and Partition Based Kronecker Products}

 The Kronecker product of two matrices $A=(a_{i,j})\in {\cal M}_{m\times n}$ and $B=(b_{r,s})\in {\cal M}_{u\times v}$ can be realized in the following two steps:
\begin{itemize}
\item[] Step 1: Produce a set of data, called the Kronecker set, as
\begin{align}\label{kd.3.2.1}
\begin{array}{l}
D_{AB}=\left\{ a_{i,j}b_{r,s}\;|\; i\in[1,m],\right.\\
~~  \left. j\in[1,n],~ r\in[1,u],~ s\in[1,v]\right\}.
\end{array}
\end{align}
 \item[] Step 2: Arrange the data in $D_{AB}$ into a matrix in the order that $Id(\a,\b;mu,nv)$, where $\vec{\a}=\vec{i}\vec{r}$ is the row index and $\vec{\b}=\vec{j}\vec{s}$ is the column index. Then the resulting matrix is $A\otimes B$.
\end{itemize}

In fact, we can define some other Kronecker products of matrices. Say,
$$
A\boxtimes B:=B\otimes A.
$$

Based on this observation, the Kronecker product of hypermatrices can also be defined through these two steps.

\begin{dfn}\label{dkd.3.2.1} Let ${\cal A}=(a_{i_1,\cdots,i_r})\in \R^{m_1\times \cdots\times m_r}$ and
${\cal B}=(b_{j_1,\cdots,j_s})\in \R^{n_1\times \cdots\times n_s}$.
\begin{itemize}
\item[(i)] The Kronecker matrix set is defined as
\begin{align}\label{kd.3.2.2}
\begin{array}{l}
D_{{\cal A}{\cal B}}=\left\{ a_{i_1,\cdots,i_r}b_{j_1,\cdots,j_s} \;|\; i_k\in[1,m_k],k\in [1,r], \right.\\
~~\left.~  j_q\in[1,n_q],~ q\in[1,s]\right\}.
\end{array}
\end{align}
\item[(ii)] A special arrange of the elements of $D_{{\cal A}{\cal B}}$ is called a Kronecker product of hypermatrices.
\end{itemize}
\end{dfn}

Let $A$ and $B$ be two matrix expressions of ${\cal A}$ and ${\cal B}$ under a pre-assigned matrix form. Then a Kronecker product of ${\cal A}$ and ${\cal B}$ can be defined by $A\otimes B$.
Next, we consider some commonly used Kronecker products of hypermatrices. They are defined by using pre-assigned special matrix forms.

\begin{itemize}
\item Outer Kronecker Product.
\end{itemize}

\begin{dfn}\label{dkd.3.2.2} \cite{lim13} Let ${\cal A}=(a_{i_1,\cdots,i_s})\in \R^{m_1\times \cdots\times m_s}$ and
 ${\cal B}=(b_{j_1,\cdots,j_t})\in \R^{n_1\times \cdots\times n_t}$. The outer-product of ${\cal A}$ and ${\cal B}$, denoted by
 \begin{align}\label{kd.3.2.3}
{\cal A}\circledcirc {\cal B}:={\cal C}=(c_{i_1,\cdots,i_s,j_1,\cdots,j_t})\in \R^{m_1\times \cdots\times m_s\times n_1\times \cdots\times n_t},
\end{align}
where
$$
c_{i_1,\cdots,i_s,j_1,\cdots,j_t}=a_{i_1,\cdots,i_s}b_{j_1,\cdots,j_t}.
$$
 \end{dfn}

 \begin{rem}\label{rkd.3.2.3}
 \begin{itemize}
\item[(i)] Although the product defined by (\ref{kd.3.2.3})  is commonly called the outer product of hypermatrices, it is easy to see that it is a Kronecker product, so we call it the outer Kronecker product.
 \item[(ii)]  The outer product of two hypermatrices is commonly denoted by $\otimes$. But to distinct it with other Kronecker product of hypermatrices, we use $\circledcirc$ instead.
 \item[(iii)] It is easy to verify that this Kronecker product is based on the matrix form $V({\cal A})$ of ${\cal A}$. That is:
\begin{align}\label{kd.3.2.4}
V({\cal A}\circledcirc {\cal B}):=V({\cal A})\otimes  V({\cal A}).
\end{align}
  \end{itemize}
 \end{rem}

\begin{itemize}
\item Partition-Based Kronecker Product.
\end{itemize}

\begin{dfn}\label{dkd.3.2.4} Let ${\cal A}=(a_{i_1,\cdots,i_r})\in \R^{m_1\times \cdots\times m_r}$ and
${\cal B}=(b_{j_1,\cdots,j_s})\in \R^{n_1\times \cdots\times n_s}$. Assume $\vec{i}=\vec{a}\vec{b}$ and $\vec{j}=\vec{c}\vec{d}$.
 Based on these two index partitions, we have spacial matrix expressions as
$$
\begin{array}{l}
A:=M^{\vec{a}\times \vec{b}}({\cal A}),\\
B:=M^{\vec{c}\times \vec{d}}({\cal B}).\\
\end{array}
$$
Then the partition-based Kronecker product of ${\cal A}$ and ${\cal B}$, denoted by
$$
{\cal C}={\cal A}\boxtimes {\cal B},
$$
is defined by
\begin{align}\label{kd.3.2.5}
M^{\vec{a}\vec{c}\times \vec{b}\vec{d}}({\cal C}):=A\otimes B.
\end{align}
\end{dfn}

Note that (\ref{kd.3.2.5}) itself shows that partition-based Kromecker product is
defined by using pre-assigned special matrix forms.

\begin{itemize}
\item Paired Kronecker Product.
\end{itemize}

 \begin{dfn}\label{dkd.3.2.5} \cite{bat16} Let ${\cal A}=(a_{i_1,\cdots,i_d})\in \R^{m_1\times \cdots\times m_d}$ and
${\cal B}=(b_{j_1,\cdots,j_d})\in \R^{n_1\times \cdots\times n_d}$.
The Kronecker product of ${\cal A}$  and ${\cal B}$, denoted by
 \begin{align}\label{kd.3.2.6}
{\cal A}\otimes_d {\cal B}:={\cal C}=(c_{k_1,\cdots,k_d})\in \R^{u_1\times\cdots\times u_d},
\end{align}
where $\vec{k}_s=\vec{i}_s\vec{j}_s$, $u_s=m_s\times n_s$,  and
 \begin{align}\label{kd.3.2.7}
c_{k_1,\cdots k_s}=a_{i_1,\cdots,i_s}b_{j_1,\cdots,j_s},\quad s\in[1,d].
\end{align}
\end{dfn}

This kind of Kronecker products is one of the most commonly accepted Kronecker products of hypermatrices.

\begin{rem}\label{rkd.3.2.6}
\begin{itemize}
 \item[(i)] In Definition \ref{dkd.3.2.5} since $\vec{k}_s=\vec{i}_s\vec{j}_s$  we have
\begin{align}\label{kd.3.2.8}
k_s=(i_s-1)n_s+j_s,\quad s\in [1,d].
\end{align}
We recall (\ref{kd.1.1.7})-(\ref{kd.1.1.8}) for the index conversion.
\item[(ii)] This definition coincides the conventional definition of Kronecker product of matrices, so it is widely accepted.
\item[(iii)] When the order of ${\cal A}$ is not the same as ${\cal B}$ the the dumy index $\emptyset$ can be used to make the orders being formally equal. In this way, the Paired Kronecker Product is applicable to two arbitrary hypermatrices.
\end{itemize}
\end{rem}

It is easy to see that this Kronecker product is independent of particular matrix expression as stated in the following proposition.

\begin{prp}\label{pkd.3.2.7}
Assume
$$
\begin{array}{l}
{\cal C}=(c_{i_1,\cdots,i_d})\in\R^{m_1n_1\times\cdots\times m_dn_d},\\
{\cal A}=(a_{r_1,\cdots,r_d})\in\R^{m_1\times\cdots\times m_d},\\
{\cal B}=(b_{s_1,\cdots,s_d})\in\R^{n_1\times\cdots\times n_d}.\\
\end{array}
$$
$$
{\cal C}={\cal A}\otimes_d {\cal B},
$$
if and only if, for any $0\leq k\leq d$
$$
\begin{array}{l}
C=M^{\vec{i}_1\cdots \vec{i}_k\times \vec{i}_{k+1}\cdots \vec{i}_d}({\cal C}),\\
A=M^{\vec{r}_1\cdots \vec{r}_k\times \vec{r}_{k+1}\cdots \vec{r}_d}({\cal A}),\\
B=M^{\vec{s}_1\cdots \vec{s}_k\times \vec{s}_{k+1}\cdots \vec{s}_d}({\cal B}),\\
\end{array}
$$
satisfy
$$
C=A\otimes B.
$$
\end{prp}

\section{Solvability of KPD of Hypermatrices}

\subsection{KPD of Outer and Partition-Based Kronecker Products}

First, we consider the outer product decomposition, which is defined as follows.

Let ${\cal A}=(a_{k_1,\cdots,k_d})\in \R^{m_1\times\cdots\times m_r\times n_1\times \cdots\times n_s}$. To find ${\cal B}=(b_{i_1,\cdots,i_r})\in \R^{m_1\times\times m_r}$ and ${\cal C}=(c_{j_1,\cdots,j_s})\in \R^{n_1\times\times n_s}$, $d=r+s$, such that
\begin{align}\label{kd.4.1.1}
{\cal A}={\cal B} \circledcirc {\cal C},
\end{align}
the index
$$
\vec{k}=\vec{i}\vec{j}.
$$
Then the following result is obvious.

\begin{prp}\label{pkd.4.1.2}
${\cal A}=(a_{k_1,\cdots,k_d})\in \R^{m_1\times\cdots\times m_r\times n_1\times \cdots\times n_s}$, $d=r+s$ is outer Kronecker product decomposable w.r.t. $r\times s$, if and only if,
$V({\cal A})$ is decomposable w.r.t. ${\bf m}\times {\bf n}$, where ${\bf m}=\prod_{k=1}^rm_k$ and ${\bf n}=\prod_{k=1}^sn_k$
\end{prp}

Now the problem is converted to the decomposition of vectors, so its solvability is  well known.

Next, we consider the partition-based KPD of hypermatrices. The following result is obvious.

\begin{prp}\label{pkd.4.1.3}
${\cal A}=(a_{k_1,\cdots,k_d})\in \R^{m_1\times\cdots\times m_d}$, where $m_k=u_kv_k$, $k\in [1,d]$,  is partition-based KPD ronecker product decomposable w.r.t. a partition
$\vec{k}=\vec{i}\vec{j}$,
if and only if,
$$
M^{\vec{i}\times \vec{j}}({\cal A})
$$
is decomposable w.r.t. $u^i\times v^i$ and $u^j\times v^j$, where
$u^i=\prod_{k\in {\bf i}}u_k$, $u^j=\prod_{k\in {\bf j}}u_k$,
$v^i=\prod_{k\in {\bf i}}v_k$, $v^j=\prod_{k\in {\bf j}}v_k$.
\end{prp}

Now the partition-based KPD of hybermatrix is converted to the matrix KPD problem, which has been completely solved before.

The paired KPD is relatively complicated. Some preparations are necessary.

\subsection{Matrix Expression of Paired Kronecker Product}

As aforementioned that  the Kronecker product of matrices can be considered as the Kronecker product of order $2$ hypermatrices under the default matrix form $A=M^{r\times n}({\cal A})$ for ${\cal A}\in \R^{m\times n}$.
Motivated by this,  it seems to be convenient to express the paired Kronecker product of hypermatrices in a ``normal form".
Then the paired KPD can also be expressed via the normal form.

To make the study easier, we first consider the case of order $3$ (i.e., cubic matrices).

Assume the normal form of ${\cal A}\in \R^{s\times m\times n}$ is as in (\ref{kd.3.1.1}). Then the following result comes from straightforward computations.

\begin{prp}\label{pkd.4.2.1}  Let $A=(a_{k,i,j})$ and $B=(b_{r,p,q})$ be the (default) matric expressions of ${\cal A}\in \R^{s\times m\times n}$ and ${\cal B}\in \R^{t\times u\times v}$ respectively. Then,
\begin{itemize}
\item[(i)]
\begin{align}\label{kd.4.2.2}
V({\cal A})=V_r(A),
\end{align}
where $V_r(A)$ is the row stacking form of $A$.

\item[(ii)] Assume ${\cal C}={\cal A}\otimes_3{\cal B}$,  $\vec{c}=\vec{k}\vec{r}$, $\vec{a}=\vec{i}\vec{p}$, $\vec{b}=\vec{j}\vec{q}$, and
$$
C=M^{\vec{c}\vec{a}\times \vec{b}}({\cal C}).
$$
Then
\begin{align}\label{kd.4.2.3}
C=\begin{bmatrix}
A^{(1)}\otimes B^{(1)}\\
A^{(1)}\otimes B^{(2)}\\
\cdots\\
A^{(1)}\otimes B^{(t)}\\
\vdots\\
A^{(2)}\otimes B^{(1)}\\
A^{(2)}\otimes B^{(2)}\\
\cdots\\
A^{(2)}\otimes B^{(t)}\\
\vdots\\
A^{(s)}\otimes B^{(1)}\\
A^{(s)}\otimes B^{(2)}\\
\cdots\\
A^{(s)}\otimes B^{(t)}\\
\end{bmatrix}.
\end{align}
\end{itemize}
\end{prp}

Motivated by the normal form of cubic matrix,  the normal form expression of a hypermatrix with order $d$ can be constructed as follows:

Let ${\cal A}=(a_{k_1,\cdots.k_d})\in \R^{n_1\times\cdots\times n_d}$. Choosing any partition $\vec{k}=\vec{i}\vec{j}$. The matrix
\begin{align}\label{kd.4.2.4}
A=M^{\vec{i}\times \vec{j}}({\cal A})\in {\cal M}_{\a\times \b},
\end{align}
where $\a=\prod_{i\in \vec{i}}n_i$ and $\b=\prod_{i\in \vec{j}}n_i$, can be used for the (default) normal  matrix form of ${\cal A}$.

Similarly to cubic matrices, setting a default form for  order $d>3$ set of hypermatrices is of first importance. Because only under a certain matrix form the paired Kronecker product can be computed by the  matrix product of their default matrix forms.  Then further manipulating hypermatrices, such as KPD etc., becomes executable.

The idea for calculating paired Kronecker product of hypermatrices via the product of their default matrix forms is as follows.
Let $A$ and $B$ be certain matrix expressions of ${\cal A}$ and ${\cal B}$ respectively. Then the matrix set of
$A\otimes B$ is the same as the matrix set of ${\cal A}\otimes {\cal B}$. Hence, we can get the matrix set of ${\cal A}\otimes {\cal B}$ from $A\otimes B$. Then we need only to arrange the elements of this matrix set in the order of the default matrix form of
${\cal A}\otimes {\cal B}$.

Assume ${\cal A}=(a_{k_1,\cdots,k_d})\in \R^{g_1\times\cdots\times g_d }$. Fix an index partition as $\vec{k}=\vec{i}\vec{j}$. Say, $\vec{i}=\vec{k}_1\cdots\vec{k}_s$ and $\vec{j}=\vec{k}_{s+1}\cdots\vec{k}_{d}$. Based on this fixed index partition we can express
\begin{align}\label{kd.4.2.5}
{\cal A}=(a_{i_1,\cdots,i_s,j_1,\cdots,j_t})\in \R^{m_1\times \cdots m_s\times n_1\times \cdots\times n_t}.
\end{align}
Similarly, another order $d$ hypmatrix ${\cal B}$ can also be expressed as
\begin{align}\label{kd.4.2.6}
{\cal B}=(b_{p_1,\cdots,p_s,q_1,\cdots,q_t})\in \R^{u_1\times \cdots u_s\times v_1\times \cdots\times v_t}.
\end{align}
Then the default matrix express of ${\cal A}$ and ${\cal B}$, based on the fixed index partition, can be expressed as
\begin{align}\label{kd.4.2.7}
\begin{array}{l}
A=M^{\vec{i}\times \vec{p}}({\cal A}),\\
B=M^{\vec{j}\times \vec{q}}({\cal B}).\\
\end{array}
\end{align}

Set
$$
{\cal C}={\cal A}\otimes_d {\cal B}.
$$
Then the default matrix form of ${\cal C}$ is
\begin{align}\label{kd.4.2.8}
C=M^{\vec{a}\times \vec{b}}({\cal C}),
\end{align}
where $\vec{a}=\vec{a}_1\cdots\vec{a}_s$, $\vec{b}=\vec{b}_1\cdots\vec{b}_t$, and $\vec{a}_{\ell}=\vec{i}_{\ell} \vec{p}_{\ell}$, $\ell\in[1,s]$,
$\vec{b}_{\ell}=\vec{j}_{\ell} \vec{q}_{\ell}$, $\ell\in[1,t]$.

Consider $C$,  its $rows$ are labeled by $\vec{i}_1\vec{p}_1\cdots \vec{i}_s\vec{p}_s$ and its columns are labeled by
$\vec{j}_1\vec{q}_1\cdots \vec{j}_t\vec{q}_t$. That is how the matrix set elements
$$
\left\{ a_{i_1\cdots i_s j_1,\cdots j_t}b_{p_1\cdots p_s q_1,\cdots q_t} \right\}
$$
are arranged in $C$.

On the other hand, we consider
$$
E:=A\otimes B,
$$
where the same matrix set elements are arranges as follows:
The rows are labeled by $\vec{i}_1\cdots\vec{i}_s\vec{p}_1\cdots\vec{p}_s$ and the columns are labeled by
 $\vec{i}_1\cdots\vec{i}_s\vec{p}_1\cdots\vec{p}_s$.

Now we know how to get the default matrix expression $C$  of the paired Kronecker product from $E$.
\begin{itemize}
\item[(i)] The row index should be changed from $\vec{i}_1\cdots\vec{i}_s\vec{p}_1\cdots\vec{p}_s$ to $\vec{i}_1\vec{p}_1\cdots \vec{i}_s\vec{p}_s$.
A permutation $\sigma\in {\bf S}_{2s}$ is required as
$$
\sigma(i)=\begin{cases}
2i-1,\quad i\in[1,s],\\
2(i-s),\quad i\in[s+1,2s].
\end{cases}
$$
Hence we construct a permutation matrix as
\begin{align}\label{kd.4.2.9}
W_L:=W^{\sigma}_{(m_1\times\cdots\times m_s\times u_1\times \cdots\times u_s)}.
\end{align}
\item[(ii)] The column index should be changed from $\vec{j}_1\cdots\vec{j}_t\vec{q}_1\cdots\vec{q}_t$ to $\vec{j}_1\vec{q}_1\cdots \vec{j}_t\vec{q}_t$.
A permutation~$\mu\in {\bf S}_{2t}$ is required as
$$
\mu(i)=\begin{cases}
2i-1,\quad i\in[1,t],\\
2(i-t),\quad i\in[t+1,2t].
\end{cases}
$$
Hence we construct a permutation matrix as
\begin{align}\label{kd.4.2.10}
W_R:=W^{\mu}_{(n_1\times\cdots\times n_t\times v_1\times \cdots\times v_t)}.
\end{align}
\end{itemize}

Then the following result is obvious.

\begin{prp}\label{pkd.4.2.2}
Consider the set of order $d$ hypermatrices, and assume an index partition
$\vec{k}=\vec{i}_1\cdots \vec{i}_s\vec{j}_1\cdots \vec{j}_t$, ($s+t=d$) is pre-assigned. The corresponding default matrix form of ${\cal A}$ is
$$
A=M^{\vec{i}\times \vec{j}}({\cal A}).
$$
Given two order $d$ hypermatrices ${\cal A}$ and ${\cal B}$ as in (\ref{kd.4.2.5}) and (\ref{kd.4.2.6}). Let
${\cal C}={\cal A}\otimes_d {\cal B}$. Then the default matrix form of ${\cal C}$ is
\begin{align}\label{kd.4.2.11}
C=W_L (A\otimes B) W^{\mathrm{T}}_R,
\end{align}
where $W_L$ and $W_R$ are defined in (\ref{kd.4.2.9}) and (\ref{kd.4.2.10}) respectively.
\end{prp}

Next, we propose a universal normal form for order $d$ hypermatrices.

 Motivated by the normal forms of $\R^{n}$, which is a column vector $V\in \R^n$; consider $\R^{m\times n}$, which is a matrix $A=(a_{i_1,\cdots,i_d})\in \R^{m\times n}$; and $\R^{s\times m\times n}$, which is a matrix $M\in {\cal M}_{sm\times n}$. We propose a universal normal form for ${\cal A}\in F^{n_1\times \cdots\times n_d}$, as follows.
$$
\vec{i}=\vec{j}\vec{k},
$$
where
$$
|\vec{j}|=
\begin{cases}
[d/2]+1,\quad d ~\mbox{is odd},\\
[d/2],\quad   d~\mbox{is even}.\\
\end{cases}
$$
Then the normal form of ${\cal A}$ is expressed as
\begin{align}\label{kd.4.2.12}
A=M^{\vec{j}\times \vec{k}}({\cal A}).
\end{align}

It is not difficult to prove that this normal form has similar properties as shown in Proposition \ref{pkd.4.2.1} for cubic matrices.
To save space, we give a simple example to describe this. Similar expression for general case can also be easily obtained.

\begin{exa}\label{ekd.4.2.3}
\begin{itemize}
\item[(i)] Consider ${\cal A}=(a_{r,s,i,j})\in \R^{m\times n\times c\times e}$. Then the normal form of ${\cal A}$ is
\begin{align}\label{kd.4.2.13}
A=\begin{bmatrix}
a^{1,1}_{1,1}&a^{1,1}_{1,2}&\cdots a^{1,1}_{1,e}\\
a^{1,1}_{2,1}&a^{1,1}_{2,2}&\cdots a^{1,1}_{2,e}\\
\cdots&~&~&~\\
a^{1,1}_{c,1}&a^{1,1}_{c,2}&\cdots a^{1,1}_{c,e}\\
a^{1,2}_{1,1}&a^{1,2}_{1,2}&\cdots a^{1,2}_{1,e}\\
a^{1,2}_{2,1}&a^{1,2}_{2,2}&\cdots a^{1,2}_{2,e}\\
\cdots&~&~&~\\
a^{1,2}_{c,1}&a^{1,2}_{c,2}&\cdots a^{1,2}_{c,e}\\
\vdots&~&~&~\\
a^{m,n}_{1,1}&a^{m,n}_{1,2}&\cdots a^{m,n}_{1,e}\\
a^{m,n}_{2,1}&a^{m,n}_{2,2}&\cdots a^{m,n}_{2,e}\\
\cdots&~&~&~\\
a^{m,n}_{c,1}&a^{m,n}_{c,2}&\cdots a^{m,n}_{c,e}\\
\end{bmatrix}=
\begin{bmatrix}
A^{(1)}\\
A^{(2)}\\
\vdots\\
A^{(mn)}
\end{bmatrix}.
\end{align}

\item[(ii)]
\begin{align}\label{kd.4.2.14}
V({\cal A})=V_r(A).
\end{align}

\item[(iii)] Assume ${\cal B}=(b_{r,s,i,j})\in \R^{p\times q\times g\times h}$ and its normal form is
$$
B=\begin{bmatrix}
B^{(1)}\\
B^{(2)}\\
\vdots\\
B^{(pq)}
\end{bmatrix}.
$$
Let ${\cal C}={\cal A}\otimes_4 {\cal B}$. Then the normal form of ${\cal C}$ is
\begin{align}\label{kd.4.2.15}
C=\begin{bmatrix}
A^{(1)}\otimes B^{(1)}\\
A^{(1)}\otimes B^{(2)}\\
\cdots\\
A^{(1)}\otimes B^{(pq)}\\
A^{(2)}\otimes B^{(1)}\\
A^{(2)}\otimes B^{(2)}\\
\cdots\\
A^{(2)}\otimes B^{(pq)}\\
\vdots\\
A^{(mn)}\otimes B^{(1)}\\
A^{(mn)}\otimes B^{(2)}\\
\cdots\\
A^{(mn)}\otimes B^{(pq)}\\
\end{bmatrix}
\end{align}
\end{itemize}
\end{exa}

\subsection{Solving Paired  KPD of Hypermatrices}

\begin{dfn}\label{dkd.4.3.1}
Given a hypermatrix ${\cal A}=(a_{k_1,\cdots,k_d})\in \R^{q_1\times \cdots\times q_d}$, where $q_k=m_kn_k$, $k\in [1,d]$.
The paired KPD problem is finding ${\cal B}=(b_{i_1,\cdots,i_d})\in \R^{m_1\times \cdots\times m_d}$ and
${\cal C}=(c_{j_1,\cdots,j_d})\in \R^{n_1\times \cdots\times n_d}$ such that
\begin{align}\label{dkd.4.3.1}
{\cal A}={\cal B}\otimes_d{\cal C}.
\end{align}
\end{dfn}

Note that both $V({\cal B}\times_d {\cal C})$ and $V({\cal B})\ltimes V({\cal C})$ have the same matrix set as
$$
\left\{b_{i_1,\cdots.i_d}c_{j_1,\cdots,j_d}\right\}.
$$
The only difference is in $V({\cal B}\times_d {\cal C})$ they are arranged in the order of
$\vec{i}_1\vec{j}_1\cdots \vec{i}_d \vec{j}_d$ while in  $V({\cal B})\ltimes V({\cal C})$ they are arranged in the order of
$\vec{i}_1\cdots \vec{i}_d\vec{j}_1\cdots\vec{j}_d$.
So a permutation $\sigma\in {\cal S}_{2d}$ is required, which is
\begin{align}\label{kd.4.3.101}
\sigma(i)=\begin{cases}
2i-1,\quad i\in[1,d],\\
2(i-d),\quad i\in[d+1,2d].
\end{cases}
\end{align}
Then the corresponding permutation matrix can be constructed as
\begin{align}\label{kd.4.3.2}
\Psi:=W^{\sigma}_{(m_1\times\cdots\times m_d\times n_1\times \cdots\times n_d)}.
\end{align}

Here $\Psi$ can convert the element arrangement in $V({\cal B})\ltimes V({\cal C}))$ to the element arrangement in
$V({\cal B}\otimes_d{\cal C})$.

Note that since $\sigma(1)=1$ and $\sigma(2d)=2d$. Then $\Psi$ can be simplified to
\begin{align}\label{kd.4.3.201}
\Psi=I_{m_1}\otimes W^{\mu}_{m_2\times\cdots m_d\times n_\times \cdots\times n_{d-1}}\otimes I_{n_d},
\end{align}
where $\mu\in {\bf S}_{2(d-1)}$ with
\begin{align}\label{kd.4.3.202}
\mu(i)=\begin{cases}
d+(i-1),\quad i ~\mbox{is odd},\\
i/2,\quad i~\mbox{is even}.
\end{cases}
\end{align}

The above argument leads to the following result.

\begin{lem}\label{lkd.4.3.2}
Assume ${\cal B}=(b_{i_1,\cdots,i_d})\in \R^{m_1\times \cdots\times m_d}$ and
${\cal C}=(c_{j_1,\cdots,j_d})\in \R^{n_1\times \cdots\times n_d}$. Then
\begin{align}\label{kd.4.3.3}
V({\cal B}\otimes_d{\cal C})=\Psi(V({\cal B})\ltimes V({\cal C})),
\end{align}
where $\Psi$ is defined in (\ref{kd.4.3.2}).
\end{lem}

Note that $\Psi$ is an orthogonal matrix, then (\ref{kd.4.3.3}) can be rewritten to
\begin{align}\label{kd.4.3.4}
\Psi^{\mathrm{T}}V({\cal B}\otimes_d{\cal C})=(V({\cal B})\ltimes V({\cal C})).
\end{align}

Now the problem of paired KPD has been converted into the corresponding KPD of vectors. Hence Lemma \ref{lkd.4.3.2} with the
result of MDA, described  in Proposition \ref{pkd.1.2.5} and the Remark \ref{rkd.1.2.501}, leads to our main result,  which is a general necessary and sufficient condition for arbitrary hypermatrix. It is so imaginary simple!

\begin{thm}\label{tkd.4.3.3} The paired KPD problem for hypermatrix ${\cal A}\in \R^{q_1\times \cdots\times q_d}$ with $q_k=m_kn_k$, $k\in [1,d]$ is solvable, if and only if, the vector
\begin{align}\label{kd.4.3.5}
V:=W^{\mathrm{T}}V({\cal A})
\end{align}
is decomposable w.r.t. $m\times n$, where $m=\ltimes_{i=1}^dm_i$ and $n= \ltimes_{i=1}^dn_i$.
\end{thm}

\begin{rem}\label{rkd.4.2.4}
\begin{itemize}
\item[(i)] Using the MDA to the vector $V$ in (\ref{kd.4.3.5}), if the $V$ is successfully decomposed into $V=h_0V_b\otimes V_c$. Then it is a straightforward to find ${\cal B}$ and ${\cal C}$ such that $V({\cal B})=V_b$ and $V({\cal C}) =V_c$.
\item[(ii)] The techniques developed for least square approximated KPD and the finite sum (precise) KPD for arbitrary  vectors can also used the corresponding problems of KPD of hypermatrices.
\item[(iii)] The paired Kronecker product can be used for two hypermatrices of different orders by adding dummy index $\emptyset$. Hence the series KPD of hypermatrices can also be done for factor hypermatrices with different orders.
 \end{itemize}
 \end{rem}

 The following corollary is easily verifiable.

 \begin{cor}\label{ckd.4.3.5}
 \begin{itemize}
 \item[(i)] If ${\cal A}={\cal B}_1\otimes_d \cdots\otimes_d {\cal B}_k$, then
\begin{align}\label{kd.4.3.6}
{\cal A}={\cal B}_1\otimes_d({\cal B}_2\otimes_d(\cdots\otimes_d {\cal B}_k)\cdot).
\end{align}
Hence, the multi-fold paired KPD  can be done step by step.
 \item[(ii)] The paired KPD is unique up to a set of constant product coefficients. That is, if ${\cal A}={\cal B}_1\otimes_d \cdots\otimes_d {\cal B}_k ={\cal C}_1\otimes_d \cdots\otimes_d {\cal C}_k$, then there is a set of coefficients $\mu_i$, $i\in [1,k]$ such that
 $$
 {\cal B}_i=\mu_i{\cal C}_i,\quad i\in [1,k],
 $$
 with $\prod_{i=1}^k\mu_i=1$.
 \end{itemize}
 \end{cor}

The following  numerical example shows the efficiency of our main result  for the paired KDP of hypermatrices.

\begin{exa}\label{ekd.4.3.6} Consider a cubic matrix ${\cal N}\in \R^{4\times 6\times 4}$ with it default matrix form
as
\begin{align}\label{kd.4.3.8}
N=(m_{k_1,k_2,k_3})=\begin{bmatrix}
N^{(1)}\\
N^{(2)}\\
N^{(3)}\\
N^{(4)}\\
\end{bmatrix}\in \R^{4\times 6\times 4},
\end{align}
where
$$
N^{(1)}=\begin{bmatrix}
6&     0&     3&     0\\
-4&    -2&    -2&    -1\\
2&     4&     1&     2\\
0&     0&     9&     0\\
0&     0&    -6&    -3\\
0&     0&     3&     6
\end{bmatrix},
$$
$$
N^{(2)}=\begin{bmatrix}
2&     0&     1&     0\\
-2&     2&    -1&     1\\
4&    -4&     2&    -2\\
0&     0&     3&     0\\
0&     0&    -3&     3\\
0&     0&     6&    -6
\end{bmatrix},
$$
$$
N^{(3)}=\begin{bmatrix}
-3&     0&     6&     0\\
2&     1&    -4&    -2\\
-1&    -2&     2&     4\\
3&     0&     3&     0\\
-2&    -1&    -2&    -1\\
1&     2&     1&     2\\
\end{bmatrix},
 $$
 $$
N^{(4)}=\begin{bmatrix}
 -1&     0&     2&     0\\
  1&    -1&    -2&     2\\
 -2&     2&     4&    -4\\
 1&     0&     1&     0\\
 -1&     1&    -1&     1\\
 2&    -2&     2&    -2\\
\end{bmatrix}.
$$

The dimension of $N$ is $4\times 6\times 4$, a guessed dimension factorization is
$$
4\times 6\times 4=2*2\times 2*3\times 2*2.
$$
(In fact, the possible factorization is not unique, one might try other factorizations.)

First we construct the permutation matrix. Assume the index
$$
\vec{k}_1\vec{k}_2\vec{k}_3=\vec{i}_1\vec{j}_1\vec{i}_2\vec{j}_2\vec{i}_3\vec{j}_2,
$$
where $i_1\in \vec{i}_1=[1,2]$, $i_2\in \vec{i}_2=[1,2]$, $i_3\in \vec{i}_3=[1,2]$,$j_1\in \vec{j}_1=[1,2]$, $j_2\in \vec{j}_2=[1,3]$, $j_3\in \vec{j}_3=[1,2]$. As discussed in the above, we need a $\sigma\in {\bf S}_6$, which permutates the order $(i_1,i_2,i_3.j_1,j_2,j_3)$
to $(i_1,j_1,i_2.j_2,i_3,j_3)$. That is,
$$
[1,2,3,4,5,6]\xrightarrow{\sigma} [1,4,2,5,3,6].
$$
Using (\ref{kd.4.3.2})-(\ref{kd.4.3.201}), we can construct $\mu$, with
$$
[1,2,3,4]\xrightarrow{\mu} [3,1,4,2].
$$
Then the permutation matrix can be constructed as
$$
\begin{array}{l}
W=I_2\otimes W^{\mu}_{[2\times 2\times 2\times 3]}\otimes I_2=\\
~~\d_{24}[1,3,5,13,15,17,24,6,14,16,18,7,\\
~~9,11,19,21,23,8,10,12,20,22,24].
\end{array}
$$
Using $V({\cal N})=V_c(N)$ and
$$
V:=W^{\mathrm{T}}V({\cal N})=(6,0,-4,-2,\cdots,-1,1,2,-2)^{\mathrm{T}}\in \R^{96}.
$$
We have $e(V)=1$ and $h_0(V)=6$. Setting $V_0=V/6$, and using MDA yield
$$
\begin{array}{l}
V_1=\Xi^e_{(1;2)}V_0=(I_8\otimes [\d^1_{12}]^{\mathrm{T}})V0=\\
(1,0.5,0,1.5,-0.5,1,0.5,0.5)^{\mathrm{T}},\\
V_2=\Xi^e_{(2;2)}V_0=([\d^1_{8}]^{\mathrm{T}}\otimes I_{12})V0=\\
~~(1,0,-0.6667,-0.3333,0.3333,0,6667,0.3333,0,\\
~~-0.3333,0.3333,0.6667,-0.6667)^{\mathrm{T}}.\\
\end{array}
$$
Finally, it is easy to verify that the following paired KPD holds:
$$
{\cal N}=6{\cal A}\otimes_3 {\cal B},
$$
where (in normal form)
$$
A=\begin{bmatrix}
1&0.5\\
0&1.5\\
-0.5&1\\
0.5&0.5
\end{bmatrix}
$$
$$
B=\begin{bmatrix}
1&0\\
-0.6667&-0.3333\\
0.3333&0,6667\\
0.3333&0\\
-0.3333&0.3333\\
0.6667&-0.6667\\
\end{bmatrix}
$$
\end{exa}

\begin{rem}\label{rkd.4.3.7} Since the problem of KPD has been converted to the KPD of vectors. If the vector $\Psi^{\mathrm{T}}V({\cal N})$ is not exact decomposable, then the least square error approximated KPD can be obtained. Say,
$$
\Psi^{\mathrm{T}}V({\cal N})\approx \ltimes_{i=1}^kV_i,
$$
where $V_i$, $i\in [1,k]$ are the least square error approximation. Then we can do the finite sum KPD by using MDA on
$$
\Psi^{\mathrm{T}}V({\cal N})- \ltimes_{i=1}^kV_i.
$$
Repeating this procedure, the precise finite sum KPD can also be obtained.

In one word, the three kinds of KPDs of arbitrary hypermatrix can be solved by using the technique developed.
\end{rem}

\section{Conclusion}

The KPD of vectors has firstly discussed in detail. The MDA is developed to solve it. When a vector is not decomposable, the MDA provides an approximate decomposition. Starting from this approximate solution, a gradient descent based algorithm is developed to provide least square error approximate decomposition. Using this algorithm again on the difference of $x$ and its approximate decomposed form several times a (precise) finite sum decomposition can be obtained. Hence the KPDs of vectors are completely solved. Then the KPD of matrix is investigated. By converting it to its equivalent KPD of the corresponding vector, the KPD is also completely solved. Finally, the matrix expression of a hypermatrix is discussed and by using its matrix expression the KPD of hypermatrix can also be converted to its equivalent KPD of corresponding vector, which also provide complete solution to KPD of hypermatrix. The results are surprisingly simple and numerically straightforward computable, comparing with existing results. We conclude that the technique developed in this paper provides a universal solution to KPD problems.

The author would like to emphasize the following fact. A hypermatrix is a matrix set but a matrix form. All the results obtained in this paper about hypermatrices are based on their matrix expressions. This is exactly the same as in Linear Algebra (or Matrix Theory), where they use particular matrix or vector form, but not the matrix set. Hence the author believes that a successful hypermatrix theory must be built on the matrix expressions of hypermatrices.

%
%

\end{document}